\documentclass[12pt, a4]{amsart}

\usepackage{latexsym}
\usepackage{amsmath}
\usepackage{amssymb}
\usepackage{amsbsy}
\usepackage[dvips]{graphicx,color}
\usepackage[all]{xy}   
\usepackage{mathrsfs}   
\def\R {\mathbb{R}}
\def\Z {\mathbb{Z}}
\def\ZZ {\mathbb{Z}/2\mathbb{Z}}
\def\C {\mathbb{C}}
\def\N {\mathbb{N}}
\providecommand{\Cls}{\mathop{\mathrm{Cls}}\nolimits}
\providecommand{\Int}{\mathop{\mathrm{Int}}\nolimits}
\title[Morse homology of nonsingular Morse-Smale foliations]{Morse homology of nonsingular Morse-Smale foliations and its duality theorem}

\author{Shigeaki Miyoshi}

\address{c/o Department of Mathematics\\ Chuo University\\
   1-13-27 Kasuga Bunkyo-ku, Tokyo\\ 112-8551, Japan}
\email{miyoshi@math.chuo-u.ac.jp}

\keywords{nonsingular Morse-Smale flow, Morse homology, Conley index, homology index braid}
\subjclass[2020]{primary~57R30, secondary~57R99, 57M99, 57R25}
\newtheorem{thm}{Theorem}[section]    
\newtheorem{lem}[thm]{Lemma}          
\newtheorem{prop}[thm]{Proposition}   
\newtheorem{cor}[thm]{Corollary}
\newtheorem{sublem}[thm]{Sublemma}
\theoremstyle{definition}
\newtheorem{defn}[thm]{Definition}     
\newtheorem*{remark}{Remark}
\newtheorem*{rem_convtn}{Remark and Convention}
\newtheorem*{note}{Note}

\newtheorem{claim}{Claim}
\newtheorem*{claim*}{Claim}
\newtheorem{case}{Case}

\begin{document}

\baselineskip=18pt

\begin{abstract}    
We define a Morse homology of the foliation which is defined by a general nonsingular Morse-Smale flow, applying the framework of the Conley index theory.  In fact, we define a chain and a boundary homomorphism, therefore, a cochain and a cohomology are obtained as well.   Furthermore, we prove a duality theorem on the homology and the cohomology.  
\end{abstract}

\maketitle

\section{Introduction} 
\label{intro}
In this paper, as a generalisation of \cite{First}, we define a homology of a nonsingular Morse-Smale foliation (NMS-foliation for short) without any assumptions.  Here, we call the foliation which consists of all orbits of a nonsingular Morse-Smale flows (NMS-flow for short) as leaves an {\em NMS-foliation}.  In \cite{First}, we imposed the assumption on the connection (or the asymptotic behavior) of the leaves between compact leaves (periodic orbits) and also non-twisting assumption around compact leaves. We will define a boundary homomorphism in the general setting, and thus, we have a homology of a general NMS-foliation.  In order to define a boundary homomorphism, the Conley index theory and its development by Conley \cite{Conley} and Franzosa \cite{Franzosa} play crucial roles.  

We have two types of $k$-chains and though those boundaries are getting confused, there are two types of homologies and also cohomologies.  Between these homologies and cohomologies, we have a duality theorem.  

This paper is organized as follows.  In Section 2 we review the Conley index theory and its development (\cite{Conley}, \cite{Franzosa}) including the definition of the Conley index, a Morse decomposition, an index filtration and a homology index braid.  We also review the fundamental facts on NMS-flows in Section 3.  In Section 4 we apply the theory reviewed in Section 2 to NMS-flows and define a boundary homomorphism in a similar way of Morse homology.  This defines a homology of an NMS-foliation in the full generality.  Section 5 is devoted to the proof of a duality theorem on our homology and cohomology of the NMS-foliation.  

\section{Morse decompositions, Conley indices and homology index braids}\label{MD_CI_HIB}
\label{Sect_2}
\subsection{Morse decompositions} 
Let $M$ be a closed $n$-manifold and $F: M\times\R\rightarrow M$ a flow on $M$.  For the simplicity, we write $x\cdot t$ instead of $F(x, t)$ for any $(x, t)\in M\times\R$. Also, we write $A\cdot T:=\{ x\cdot t\ |\ x\in A, t\in T\}$ for $A\subset M$ and $T\subset\R$.  

\begin{defn}
A subset of $M$ is said to be {\em invariant} if it is a union of orbits of $F$.  
For any subset $N\subset M$, we define the {\em maximal invariant set} $I(N)$ of $N$ for $F$ to be 
\[
I(N) := \{ x\in N\ |\ x\cdot\R\subset N\}. 
\]
Note that $I(N) = \bigcap_{t\in\R}N\cdot t$.  
\end{defn}

\begin{defn}
A compact invariant subset $S\subset M$ is said to be {\em isolated} if there exists a compact neighbourhood $N$ of $S$ such that $I(N) = S$.  In the case, $N$ is called an {\em isolating neighbourhood} for $S$.  
\end{defn}

We define the {\em backward flow} $F^{\ast}$ of $F$ by $F^{\ast}(x, t):=F(x, -t)$.  

\begin{defn}
 For a subset $A\subset M$, we define the {\em $\omega$-limit set} $\omega (A)$ to be the maximal invariant set of $\Cls(A\cdot [0, +\infty))$:
 $\omega (A)=I(\mathrm{Cls}(A\cdot [0, +\infty)))$, 
where $\mathrm{Cls}$ denotes the closure.  Similarly the {\em $\alpha$-limit set} $\alpha (A)$ is defined to be the $\omega$-limit set of $A$ with respect to the backward flow $F^{\ast}$, i.e., $\alpha (A)=I(\mathrm{Cls}(A\cdot (-\infty, 0])$ with respect to $F$.  
\end{defn}

\begin{remark}
\begin{enumerate}
 \item Clearly $\bigcup_{x\in A}\omega (x)\subset\omega (A)$, however, $\bigcup_{x\in A}\omega (x)\ne\omega (A)$ in general.  
 \item It can be easily seen that $\omega (A)\left( =\bigcap_{t\in\R}\Cls(A\cdot [0, \infty))\cdot t\right) =\bigcap_{t\geq 0}\Cls(A\cdot [t, \infty))$.  
\end{enumerate}
\end{remark}

\begin{defn}
 Suppose that $S_{1}$ and $S_{2}$ are compact invariant subsets.  Let $C(S_{1}, S_{2})$ denote the subset of all points $x\in M$ whose $\alpha$-limit set $\alpha (x)$ is contained in $S_{1}$ and the $\omega$-limit set $\omega (x)$ is contained in $S_{2}$.  We call $C(S_{1}, S_{2})$ the {\em set of orbits connecting $S_{1}$ and ${S_{2}}$}.  
\end{defn}

Let $\mathcal{M} = \{ m(\pi)\}_{\pi\in P}$ be a family of mutually disjoint isolated compact invariant subsets of $(M, F)$.  Here, $P$ is an index set and we suppose that a partial order $\leq$ on $P$ is given.  In the case that $F$ is an NMS-flow, we have a natural partial order on $P$ defined by the flow (see Section \ref{def_bdry}).  

\begin{defn}
A family $\mathcal{M}=\{ m(\pi)\}_{\pi\in P}$ as above is a {\em Morse decomposition} for $(M, F)$ if for any $x\in M\setminus\cup\mathcal{M}$, there exist $\pi, \pi'\in P$ such that $\pi' <\pi$ and $x\in C(m(\pi), m(\pi'))$.  
\end{defn}
\noindent
In this paper, we assume that $C(m(\pi), m(\pi'))\ne\emptyset$ if and only if $\pi'\leq\pi$.  In the case of NMS-flow, this condition is satisfied.  

\subsection{Intervals in the index set $P$} 
In this subsection, we show some simple notions of a poset which we will use in order to define the  boundary homomorphisms of our homology.  

\begin{defn}
\begin{enumerate}
 \item A subset $I\subset P$ is an {\em interval} if for $\pi, \pi'\in I$ and $\pi''\in P$ with $\pi'<\pi''<\pi$ it follows $\pi''\in I$.  The set of all intervals in $P$ is denoted by $\mathscr{I}$.  
 \item A subset $I$ in $P$ is an {\em attracting interval} if for $\pi\in I$ and $\pi'\in P$ with $\pi'<\pi$ it follows $\pi'\in I$. The set of all attracting intervals in $P$ is denoted by $\mathscr{A}$.  
\end{enumerate}
\end{defn}
\begin{prop}
 \begin{enumerate}
  \item $\mathscr{A}\subset\mathscr{I}$. 
  \item $\emptyset, P\in\mathscr{A}$. 
  \item If $I, J\in\mathscr{I}$, then $I\cap J\in\mathscr{I}$. 
  \item If $I, J\in\mathscr{A}$, then $I\cap J, I\cup J\in\mathscr{A}$. 
 \end{enumerate}
\end{prop}

\begin{defn} 
 An {\em adjacent $m$-tuple of intervals} $(I_{1}, I_{2}, \ldots, I_{m})$ in $P$ is an $m$-tuple of mutually disjoint subsets of $P$ which satisfies the following:
\begin{enumerate}
 \item $\cup_{j=1}^{m}I_{j}\in\mathscr{I}$. 
 \item If $\pi_{j}\in I_{j}, \pi_{k}\in I_{k}$ and $j<k$, then $\pi_{k}\nless\pi_{j}$. 
\end{enumerate}
\noindent
The set of all adjacent $m$-tuples is denoted by $\mathscr{I}_{m}$. 
\end{defn}

\begin{rem_convtn}
 \begin{enumerate}
  \item $\mathscr{I}_{1}=\mathscr{I}$. 
  \item Note that in the definition of adjacent $m$-tuple, each $I_{p}$ is not assumed to be an interval.  However, it follows it is an interval automatically (see Proposition \ref{sub_adjacent_tuple}).  
  \item For any $(I_{1}, I_{2}, \ldots, I_{m})\in\mathscr{I}_{m}$, we write $I_{1}I_{2}\cdots I_{m}:=\cup_{j=1}^{m}I_{j}$. 
 \end{enumerate}
\end{rem_convtn}

\begin{defn}
 \begin{enumerate}
  \item If the both of $(I, J)$ and $(J, I)$ are adjacent pair, then we say that $I$ and $J$ are {\em non-comparable}. 
  \item If $(I_{1}, I_{2}, \ldots, I_{m})\in\mathscr{I}_{m}$, then it is said to be a {\em decomposition} of the interval $I=I_{1}I_{2}\cdots I_{m}$.   
 \end{enumerate}
\end{defn}

\begin{remark}
 If intervals $I$ and $J$ are non-comparable, then for any $\pi_{I}\in I$ and $\pi_{J}\in J$ we have $\pi_{J}\nless\pi_{I}$ since $(I, J)\in\mathscr{I}_{2}$ and $\pi_{I}\nless\pi_{J}$ since $(J, I)\in\mathscr{I}_{2}$.  That is, $\pi_{I}$ and $\pi_{J}$ are non-comparable. 
\end{remark}

\begin{prop}\label{sub_adjacent_tuple}
 Let $(I_{1}, I_{2}, \ldots, I_{m})\in\mathscr{I}_{m}$ and $p, q\in\{ 1, 2, \ldots, m\}$ with $p\leq q$.  Then we have $I_{p}I_{p+1}\cdots I_{q}\in\mathscr{I}$.  In particular, $I_{r}\in\mathscr{I}$ for any $r\in\{ 1, 2, \ldots m\}$.  
\end{prop}

The construction in the proof of the following is one of the keys in this paper:

\begin{lem}\label{K=IJ}
 For any $J\in\mathscr{I}$ there exists $I\in\mathscr{A}$ such that $IJ\in\mathscr{A}$ and $(I, J)$ is a decomposition of $IJ$. 
\end{lem}
\begin{proof}
 Let $I$ be the subset of $P\setminus J$ each element of which has an element of $J$ greater than it.  That is, $\pi\in I$ if and only if $\pi\not\in J$ and there is $\pi'\in J$ such that $\pi<\pi'$.  Then it is clear that $I, IJ\in\mathscr{A}$.  To show that $(I, J)\in\mathscr{I}_{2}$, suppose $\pi_{I}\in I$ and $\pi_{J}\in J$ with $\pi_{J}< \pi_I$ are given.  By the definition of $I$ there exists $\pi\in J$ such that $\pi_{I}<\pi$.  Then we have $\pi_{J}< \pi_{I}< \pi$ which implies $\pi_{I}\in J$ since $J\in\mathscr{I}$. This is a contradiction to the definition of $I$.     
\end{proof}

\begin{prop}\label{condensed_tuple}
 Let $(I_{1}, I_{2}, \ldots , I_{m})\in\mathscr{I}_{m}$ and $p, q\in\{ 1, 2, \ldots , m\}$ with $p<q$.  Then, we have the following;
\begin{enumerate}
 \item $(I_{1}, \ldots , I_{p-1}, \cup_{j=p}^{q}I_{j}, I_{q+1}, \ldots , I_{m})\in\mathscr{I}_{m-q+p}$.  
\item $(I_{p}, \ldots , I_{q})\in\mathscr{I}_{q-p+1}$
\end{enumerate}
\end{prop}

\begin{remark}
 By Proposition \ref{condensed_tuple}, if $(I, J, K)$ is an adjacent triple, then $(I, J), (J, K), (IJ, K), (I, JK)$ are all adjacent pairs.  
\end{remark}

\subsection{Morse sets and Morse decompositions} 
We show a summary of the Conley index theory and one of its developments due to Conley \cite{Conley} and Franzosa \cite{Franzosa} in the rest of this section.  Let $\mathcal{M}=\{ m(\pi)\}_{\pi\in P}$ be a Morse decomposition for $(M, F)$.  

\begin{defn}
 An invariant set $A$ is an {\em attractor} if there exists a neighbourhood $U$ such that $\omega (U)=A$.  Similarly $A$ is a {\em repeller} if $\alpha (U)=A$ for a neighbourhood $U$.  
\end{defn}

\begin{remark}
 The attractor we define as above might be an attractor in a weak sense.  Here, we do not impose its transitivity or minimality.  In a sense, the definition is necessary and sufficient for our purpose.  See Section \ref{NMS_RHD}, \ref{def_bdry}.  
\end{remark}

\begin{defn}
 For each $I\in\mathscr{I}$, we define 
\[
 m(I):=\left( \bigcup_{\pi\in I}m(\pi)\right)\cup\left(\bigcup_{\pi, \pi'\in I}C(m(\pi), m(\pi'))\right).  
\]
We call $m(I)$ the {\em Morse set} of $I$.  
\end{defn}

Now, suppose that $S$ is a compact invariant set and $A\subset S$ is an attractor in $S$.  Set $A^{\ast}:=\{ x\in S\ |\ \omega(x)\cap A=\emptyset\}$ and $D(A):=\{ x\in S\ |\ \omega (x)\subset A\}$.  

\begin{prop}\label{A-R_pair}
 We have the following:
\begin{enumerate}
 \item $S\setminus D(A)=A^{\ast}$,
 \item\label{cmplmnt_repeller} $A^{\ast}$ is a repeller in $S$, 
 \item\label{D(A)} $D(A)\setminus A=C(A^{\ast}, A)$. 
\end{enumerate}
\end{prop}
\begin{proof}
First, for any given $x\in S$, we show $\omega (x)\cap A\ne\emptyset$ if and only if $\omega (x)\subset A$.  Let $U$ be a neighbourhood of $A$ such that $\omega (U)=A$ and assume $\omega (x)\cap A\ne\emptyset$, then there exists $t\in\R$ such that $x\cdot t\in U$.  Thus $\omega (x)=\omega (x\cdot t)\subset\omega (U)=A$.  The converse is clear.  Hence for any $x\in S$, $x\in A^{\ast}$ or $x\in D(A)$ exclusively, which implies $S=A^{\ast}\sqcup D(A)$.  

To prove (\ref{cmplmnt_repeller}) we claim the following:

\setcounter{claim}{0}
\begin{claim}\label{1st_claim}
 There exists a positive real number $\bar{t}$ such that $\Cls(U\cdot [\bar{t}, \infty ))\subset\Int(U)$, where $\Int$ denotes the interior.  
\end{claim}
\begin{proof}
 Otherwise, for any $t>0$ there exists an $x_{t}\in (S\setminus\Int(U))\cap\Cls(U\cdot [t, \infty ))$.  Hence we have a sequence $\{ x_{n}\}_{n\in\N}$ in $S\setminus\Int(U)$ with $x_{n}\in\Cls(U\cdot [n, \infty ))$.  Since $S\setminus\Int(U)$ is compact, we have a convergent subsequence $\{ x_{n(k)}\}_{k\in\N}$ of $\{ x_{n}\}$.  Denote $x_{0}:=\lim_{k\rightarrow\infty}x_{n(k)}$.  Then $x_{0}\in S\setminus\Int(U)$.  However, we will show it lies in $A$ as well, which is a contradiction.  To show $x_{0}\in A$ note that $A=\omega (U)=\bigcap_{t\geq 0}\Cls(U\cdot [t, \infty ))$ and $\{\Cls(U\cdot [n(k), \infty))\}_{k\in\N}$ is a decreasing sequence of subsets and that for any neighbourhood $N$ of $x_{0}$, there is $\ell\in\N$ such that $x_{n(m)}\in N$ if $m\geq \ell$.  Hence, if $m\geq \ell$ then $x_{n(m)}\in N\cap\Cls(U\cdot [n(m), \infty ))$.  Therefore, for any $t>0$, $N\cap\Cls(U\cdot [t, \infty ))\ne\emptyset$, which implies $x_{0}\in\Cls(U\cdot [t, \infty ))$ for any $t>0$.  Thus, $x_{0}\in\bigcap_{t\geq 0}\Cls(U\cdot [t, \infty ))=\omega (U)=A$.  
\end{proof}
\noindent
We set $V^{\ast}:=S\setminus U\cdot [\bar{t}, \infty )$.  
\begin{claim}\label{2nd_claim}
 $\Cls(S\setminus U)\subset\Int(V^{\ast})$.  
\end{claim}
\begin{proof}
 Since $\Cls(U\cdot [\bar{t}, \infty))\subset\Int(U)$ we have
\begin{eqnarray*}
 \Cls(S\setminus U)
 & = & S\setminus\mathrm{Int (U)} \\
 & \subset & S\setminus\Cls(U\cdot [\bar{t}, \infty)) \\
 & = & \Int(S\setminus U\cdot [\bar{t}, \infty)) \\
 & = & \Int(V^{\ast})
\end{eqnarray*}
\end{proof}
\noindent
Then, for any $x\in V^{\ast}=S\setminus U\cdot [\bar{t}, \infty )$ and any $s\geq\bar{t}$, we have $x\cdot (-s)\notin U\cdot [\bar{t}-s, \infty )$, in particular, $x\cdot (-s)\notin U$.  Thus we have $V^{\ast}\cdot (-\infty, -\bar{t}]\subset S\setminus U$ and then $\Cls(V^{\ast}\cdot (-\infty, -\bar{t}])\subset\Cls(S\setminus U)$.  By Claim \ref{2nd_claim}, we have $\Cls(V^{\ast}\cdot (-\infty, -\bar{t}])\subset\Int(V^{\ast})$, therefore, $\alpha (V^{\ast})\subset\Int(V^{\ast})$.  This implies that $\alpha (V^{\ast})$ is a repeller.  

To complete the proof of (\ref{cmplmnt_repeller}) we have to show that $\alpha (V^{\ast})=A^{\ast}$.  For any $x\in\alpha (V^{\ast})$ we have $\omega (x)\subset\alpha (V^{\ast})$ since $\alpha (V^{\ast})$ is compact invariant.  Hence $\omega (x)\cap A=\emptyset$ which implies $x\in A^{\ast}$.  Thus $\alpha (V^{\ast})\subset A^{\ast}$.  Finally, for any $x\in A^{\ast}$, since $x\cdot\R\cap U=\emptyset$, we have $x\cdot\R\subset V^{\ast}$.  This implies $x\in I(\Cls(V^{\ast}\cdot (-\infty, 0]))=\alpha (V^{\ast})$.  

As for (\ref{D(A)}), since $A$ is invariant, $C(A^{\ast}, A)\subset D(A)\setminus A$ is clear.  If $x\in D(A)\setminus A$, then $\omega (x)\subset A$ and $\alpha (x)\not\subset A$, i.e., $\alpha (x)\cap A=\emptyset$.  We claim here that $\alpha (x)\subset A^{\ast}$.  Otherwise, there exists $y\in\alpha (x)-A^{\ast}$, then $y\in D(A)$, hence $\omega (y)\subset A$.  However, $\omega (y)\subset\alpha (x)$ since $y\in\alpha (x)$ and $\alpha (x)$ is compact invariant.  Therefore, we have $\omega (y)\subset\alpha (x)$ and $\omega (y)\subset A$, a contradiction to $\alpha (x)\cap A=\emptyset$, which shows $\alpha (x)\subset A^{\ast}$.  This completes the proof of Proposition \ref{A-R_pair}.  
\end{proof}

\begin{defn}
 The pair $(A, A^{\ast})$ as in Proposition \ref{A-R_pair} is called an {\em attractor-repeller pair} in $S$ and $A^{\ast}$ is called the {\em complementary repeller} to $A$.  
\end{defn}

The following is one of the crucial facts in this paper.

\begin{prop}\label{AI=>A}
 Let $S$ be an isolated compact invariant set and $\mathcal{M}=\{ m(\pi)\ |\ \pi\in P\}$ a Morse decomposition of $S$.  If $I$ is an attracting interval in $P$, then $m(I)$ is an attractor in $S$.  
\end{prop}
\noindent
\textit{Proof of Proposition \ref{AI=>A}.}
 We prove the proposition by the induction on the number of elements of $P$, $\# P$.  In the case that $P$ is a singleton, $P=\{ \theta\}$, then $m(\theta)=S$ which is an attractor in $S$.  Assume that the assertion holds if $\# P=m-1$ and now $\# P=m$.  For a given attracting interval $I$ in $P$, take $\theta\in I$ which is a local minimum in $P$, i.e., there is no $\pi\in P$ such that $\pi<\theta$.  

\begin{claim*}
 $m(\theta)$ is an attractor in $S$.  
\end{claim*}
\noindent
\textit{Proof.}
 Take an isolating neighbourhood $U$ of $m(\theta)$.  By definition, $U$ is a compact neighbourhood of $m(\theta)$ and $I(U)=m(\theta)$.  We can take $U$ so that $U\cap (\cup_{\pi\in P\setminus\{\theta\} }m(\pi))=\emptyset$.  Suppose any $x\in U\setminus m(\theta)$ is given.  Since $\mathcal{M}$ is a Morse decomposition and $x\not\in\cup_{\pi\in P}m(\pi)$, there exist $\pi, \pi'\in P$ such that $\pi' <\pi$ and $x\in C(m(\pi), m(\pi'))$.  Then $\alpha (x)\subset m(\pi)$ and $\theta\ne\pi$ since $\theta$ is a local minimum.  Thus we have $\alpha (x)\cap U=\emptyset$, which implies that $m(\theta)=I(U)$ is an attractor by the following lemma.  

\begin{lem}[\cite{Conley} II 5.1 D]
 Let $U$ be a compact subset of $S$ with the property that $\alpha (x)\cap U=\emptyset$ for any $x\in\mathrm{Fr}(U)$, where $\mathrm{Fr}$ denotes the frontier.  Then $I(U)$ is an attractor.  
\end{lem}
\hfill\textit{QED of Claim.}

\noindent
Then it can be easily seen that $m(P\setminus\{\theta\})$ is the complementary repeller to $m(\theta)$ in $S$.  Now, $\{ m(\pi)\ |\ \pi\in P\setminus\{\theta\}\}$ is a Morse decomposition of $m(P\setminus\{\theta\})$ (see Proposition \ref{MorseDecomp_MS}) and $I\setminus\{\theta\}$ is an attracting interval in $P\setminus\{\theta\}$.  Since $\# (P\setminus\{\theta\} )=m-1$, $m(I\setminus\{\theta\})$ is an attractor in $m(P\setminus\{\theta\})$ by the induction hypothesis.  Then the complementary repeller to $m(I\setminus\{\theta\})$ in $m(P\setminus\{\theta\})$ is $m(P\setminus I)$, consequently it follows that $m(P\setminus\{\theta\})$ is a repeller in $S$ and $m(P\setminus I)$ is a repeller in $m(P\setminus\{\theta\})$.  Hence, $m(P\setminus I)$ is a repeller in $S$ by the following lemma.  
\begin{lem}[\cite{Conley} II 5.3 D]
 If a subset $A$ of $S$ is an attractor $($resp. a repeller$)$ in $S$ and a subset $B$ of $A$ is an attractor $($resp. a repeller$)$ in $A$, then $B$ is an attractor $($resp. a repeller$)$ in $S$.  
\end{lem}
\noindent
Then $m(I)$ is the complementary attractor to $m(P\setminus I)$ in $S$.  In particular, $m(I)$ is an attractor in $S$.  This completes the proof.
\hfill\textit{QED of Proposition \ref{AI=>A}.}

\vspace{3ex}
Every Morse set is an intersection of an attractor and a repeller.  In fact, we have the following.  

\begin{prop}\label{MS_inv}
 For any $J\in\mathscr{I}$, there exist $K, I\in\mathscr{A}$ such that $m(J)=m(K)\cap m(I)^{\ast}$, where $m(I)^{\ast}$ is the complementary repeller to $m(I)$ in $S$.  In particular, every Morse set is isolated compact invariant set.  
\end{prop}
\begin{proof}
 First, by Lemma \ref{K=IJ} there exist $K, I\in\mathscr{A}$ such that $(I, J)$ is a decomposition of $K$.  Then, since $m(I)^{\ast}=\{ x\in S\ |\ \omega (x)\cap m(I)=\emptyset\}$ and 
\begin{eqnarray*}
 m(K) & = & \left( \bigcup_{\pi\in K}m(\pi)\right)\cup\left(\bigcup_{\pi, \pi'\in K}C(m(\pi), m(\pi')) \right)
\end{eqnarray*}
we have 
\begin{eqnarray*}
 m(K)\cap m(I)^{\ast} & = & \left( \bigcup_{\pi\in K\setminus I}m(\pi)\right)\cup\left(\bigcup_{\pi\in K, \pi'\in K\setminus I}C(m(\pi), m(\pi'))\right) \\
   & = & m(J)
\end{eqnarray*}
Note that $C(m(\pi), m(\pi'))=\emptyset$ for $\pi\in I$ and $\pi'\in J=K\setminus I$.  
Since attractors and repellers are isolated compact invariant sets, any Morse set is isolated compact invariant as well.  
\end{proof}

The following proposition is clear.  
\begin{prop}\label{MorseDecomp_MS}
 For any interval $I\in\mathscr{I}$, we have the following:
\begin{enumerate}
 \item $\{ m(\pi )\}_{\pi\in I}$ is a Morse decomposition of $m(I)$,
 \item $\{ m(\pi )\ |\ \pi\in P\setminus I\}\cup\{ m(I)\}$ is a Morse decomposition of $S$ with its index set \\ $(P\setminus I)\cup\{ I\}$ whose order is naturally induced from $P$.  
\end{enumerate} 
\end{prop}

As a consequence of Proposition \ref{MS_inv}, \ref{MorseDecomp_MS}, we can restrict Morse decompositions to Morse sets, i.e., we can only concern Morse sets in order to consider a Morse decomposition.  

\begin{cor}\label{I2=>AP}
 For any $(I, J)\in\mathscr{I}_{2}$, $(m(I), m(J))$ is an attractor-repeller pair in $m(IJ)$.  
\end{cor}
\begin{proof}
 By Proposition \ref{MorseDecomp_MS}, It follows that $\{ m(I), m(J)\}$ is a Morse decomposition of $m(IJ)$ with its index set $\{ I, J\}$.  Since $\{ I\}$ is an attracting interval in $\{ I, J\}$, $m(I)$ is an attractor in $m(IJ)$ by Proposition \ref{AI=>A}.  Then
\begin{eqnarray*}
 m(I)^{\ast} & = & \{ x\in m(IJ)\ |\ \omega (x)\cap m(I)=\emptyset\} \\
             & = & m(IJ)\setminus (m(I)\cup C(m(J), m(I))) \\
             & = & m(J).  
\end{eqnarray*}
\end{proof}

\subsection{Conley indices and index filtrations} 
\begin{defn}
An {\em index pair} $(N, L)$ for an isolated compact invariant set $S$ is a pair of compact subsets $L\subset N$ such that the following conditions are satisfied:
\begin{enumerate}
\item \label{IP_plcmnt}
$S=I(\Cls(N\setminus L))\subset\Int(N\setminus L)$.  
\item\label{pos_inv}
If $x\in L$ and $x\cdot [0, t]\subset N$, then $x\cdot [0, t]\subset L$.  
\item\label{exit_set}
If $x\in N\setminus L$, then there is a positive number $t$ such that $x\cdot [0, t]\subset N$.
\end{enumerate}
The condition (\ref{pos_inv}) implies $L$ is positively invariant relative to $N$ by the flow and the condition (\ref{exit_set}) implies that each orbit that leaves $N$ must go through $L$ first.  $L$ is called the {\em exit set} of the index pair.  Also, the quotient pointed space $N/L$ is called an {\em index space} of $S$.  
\end{defn}

\begin{remark}
 Note that the condition (\ref{exit_set}) is equivalent under the condition (\ref{IP_plcmnt}), (\ref{pos_inv}) to the following condition:\\
\vspace{-3ex}
\begin{enumerate}
 \item [$(3')$]$x\in N$ and $x\cdot\R_{\geq 0}\not\subset N$ imply that there is a non-negative real number $t$ such that $x\cdot [0, t]\subset N$ and $x\cdot t\in L$.  
\end{enumerate}
\end{remark}

The following results were proven by Conley in \cite{Conley} and by Salamon in \cite{Salamon} with new and simplified definitions and proofs.  

\begin{prop}\label{existence_indexpair}
Every isolated compact invariant set admits an index pair. 
\end{prop}

\begin{prop}\label{homotopy_inv_indexpair}
If $S$ is an isolated compact invariant set and $(N, L)$ and $(N', L')$ are two index pairs for $S$, then $N/L$ and $N'/L'$ are homotopy equivalent as pointed spaces via maps that are induced by the flow.  
\end{prop}

\begin{defn}
The {\em Conley index} of an isolated compact invariant set $S$ is the homotopy type of $N/L$ where $(N, L)$ is an index pair for $S$.  
\end{defn}
\begin{defn}
 Let $S$ be an isolated compact invariant subset and $(N, L)$ an index pair for $S$. The {\em homology index} of $S$ denoted by $\mathscr{H}_{\ast}(S)$ is defined to be $H_{\ast}(N/L)$.  
\end{defn}

We say that an index pair $(N, L)$ is {\em regular} if the inclusion map $L\hookrightarrow N$ is a cofibration.  It is known that any index pair $(N, L)$ can be modified into a regular index pair $(N, \tilde{L})$ (see Section 5.1 of \cite{Salamon}).  However, in the case of NMS-flow which is our main concern, we need not take care of the regularity.  In fact, all the minimal isolated compact invariant sets appearing in our setting are periodic orbits with each of which a {\em round handle} is associated and the round handle is a regular index pair for the periodic orbit (see Section \ref{NMS_RHD}).  Suppose that $S$ is an isolated compact invariant set and $(N, L)$ is a regular index pair for $S$.  Since for any cofibration $B\hookrightarrow Y$ with a closed subset $B$, the quotient map $(Y, B)\rightarrow (Y/B, \ast)$ induces isomorphisms between homology (resp. cohomology) groups of $(Y, B)$ and $(Y/B, \ast)$,  $H_{\ast}(N, L)$ (resp. $H^{\ast}(N, L)$) are topological invariants for the isolated compact invariant set $S$.  We refer the reader to \cite{Conley}, \cite{Salamon} and \cite{Banyaga-Hurtubise} for the proofs of Propositions \ref{existence_indexpair} and \ref{homotopy_inv_indexpair}.  

\vspace{2ex}
The following is another key lemma in this paper.
\begin{lem}\label{index_triple}
 Let $S$ be an isolated compact invariant set in $M$ and $(A, A^{\ast})$ an attractor-repeller pair in $S$.  Suppose that $N_{0}\subset N_{1}\subset N_{2}$ are compact subsets of $M$ and that $(N_{1}, N_{0})$ is an index pair for $A$ and $(N_{2}, N_{0})$ is an index pair for $S$.  Then, $(N_{2}, N_{1})$ is an index pair for $A^{\ast}$.  
\end{lem}
\noindent
\textit{Proof.}\\
\textbf{Condition (3)} 
Suppose that $x\in N_{2}\setminus N_{1}$ is given.  Then $x\in N_{2}\setminus N_{0}$ and since $(N_{2}, N_{0})$ is an index pair for $S$, by the condition (3) of $(N_{2}, N_{0})$, there exists a $t>0$ such that $x\cdot [0, t]\subset N_{2}$.  \\
\textbf{Condition (2)}
Let $x\in N_{1}$ and $x\cdot [0, t]\subset N_{2}$.  If $x\in N_{0}$ then since $(N_{2}, N_{0})$ is an index pair for $S$, by the condition (2) of $(N_{2}, N_{0})$, $x\cdot [0, t]\subset N_{0}\subset N_{1}$.  
If $x\not\in N_{0}$ then we have $x\in N_{1}\setminus N_{0}$, and by the condition (3) of the index pair $(N_{1}, N_{0})$ for $A$, there exists a $s>0$ such that $x\cdot [0, s]\subset N_{1}$.  Set $\bar{s}:=\mathrm{sup}\{ s>0\ |\ x\cdot [0, s]\subset N_{1}\}$.  Then since $N_{1}$ is compact, we have $x\cdot\bar{s}\in N_{1}$.  Now we only have to show the following.  
\setcounter{claim}{0}
\begin{claim}\label{t=<sup}
 $t\leq \bar{s}$.  
\end{claim}
\noindent
\textit{Proof of Claim \ref{t=<sup}.} Assume $t>\bar{s}$.  If $x\cdot\bar{s}\in N_{0}$, then since $(x\cdot\bar{s})\cdot [0, t-\bar{s}]=x\cdot [\bar{s}, t]\subset N_{2}$, by the condition (2) of the index pair $(N_{2}, N_{0})$ for $S$, $(x\cdot\bar{s})\cdot [0, t-\bar{s}]\subset N_{0}\subset N_{1}$.  This contradicts the definition of $\bar{s}$.  If $x\cdot\bar{s}\not\in N_{0}$, then we have $x\cdot\bar{s}\in N_{1}\setminus N_{0}$ and by the condition (3) of the index pair $(N_{1}, N_{0})$ for $A$, there exists $u>0$ such that $(x\cdot\bar{s})\cdot [0, u]\subset N_{1}$ which contradicts the definition of $\bar{s}$.  Hence we have $t\leq\bar{s}$.  
\hfill\textit{QED of Claim \ref{t=<sup}.}\\
\textbf{Condition (1)}
First, we prove that $A^{\ast}\subset\Int(N_{2}\setminus N_{1})$.  
\begin{claim}\label{A*N1}
 $A^{\ast}\cap N_{1}=\emptyset$.  
\end{claim}
\noindent
\textit{Proof.} Otherwise, there exists $x\in A^{\ast}\cap N_{1}$.  Then we have $x\cdot\R_{\geq 0}\subset A^{\ast}\subset S\subset N_{2}$ and since we have already proven that $N_{1}$ is positively invariant relative to $N_{2}$, this implies $x\cdot\R_{\geq 0}\subset N_{1}$.  Thus $x\cdot\R_{\geq 0}\subset A^{\ast}\cap N_{1}$, therefore $\omega (x)\subset A^{\ast}\cap N_{1}$.  Since $A^{\ast}\subset S\subset\Int(N_{2}\setminus N_{0})(=\Int(N_{2})-N_{0})$, we have $A^{\ast}\cap N_{0}=\emptyset$.  Hence $\omega (x)\cap N_{0}=\emptyset$, and since $\omega (x)$ is closed, we have $\omega (x)\subset\Cls(N_{1}\setminus N_{0})$.  Thus we have $\omega (x)\subset I(\Cls(N_{1}\setminus N_{0}))=A$, however, this contradicts the assumption $x\in A^{\ast}$.  

\hfill\textit{QED of Claim \ref{A*N1}.}  

\vspace{2ex}
\noindent
Since $A^{\ast}\subset S\subset\Int(N_{2}\setminus N_{0})$, Claim \ref{A*N1} implies $A^{\ast}\subset\Int(N_{2}\setminus N_{1})$.  

Finally, we prove that $A^{\ast}=I(\Cls(N_{2}\setminus N_{1}))$.  We need the following sublemma:

\begin{sublem}\label{Nested_IsoNbd}
 Let $N$ be an isolating neighbourhood of $S$ and $N'$ a compact neighbourhood of $A^{\ast}$ with $N'\cap A=\emptyset$ and $N'\subset N$.  Then $N'$ is an isolating neighbourhood of $A^{\ast}$.  
\end{sublem}
\noindent
\textit{Proof.} We have to show that $A^{\ast}=I(N')$.  First, since $A^{\ast}$ is invariant and $A^{\ast}\subset N'$, we have $A^{\ast}\subset I(N')$.  On the converse, suppose $x\in S\setminus A^{\ast}$ is given.  Recall that $S\setminus A^{\ast}=A\sqcup C(A^{\ast}, A)$.  If $x\in A$, then $x\not\in N'$ since $N'\cap A=\emptyset$.  If $x\in C(A^{\ast}, A)$, then $\omega (x)\subset A$, therefore $\omega (x)\cap N'=\emptyset$, thus $x\cdot\R\not\subset N'$ which implies $x\not\in I(N')$.  Consequently, we have $S\setminus A^{\ast}\subset S\setminus I(N')$, i.e., $A^{\ast}\supset I(N')$, which completes the proof.  
\hfill\textit{QED of Sublemma \ref{Nested_IsoNbd}.}  

\vspace{2ex}
\noindent
Now, since $A\subset\Int(N_{1})$, we have $A\cap\Cls(N_{2}\setminus N_{1})=\emptyset$.  Also, as we have shown above, $A^{\ast}\subset\Int(N_{2}\setminus N_{1})$.  Thus $\Cls(N_{2}\setminus N_{1})$ is a compact neighbourhood of $A^{\ast}$.  Clearly, $\Cls(N_{2}\setminus N_{1})\subset\Cls(N_{2}\setminus N_{0})$ and $\Cls(N_{2}\setminus N_{0})$ is an isolating neighbourhood of $S$ since $(N_{2}, N_{0})$ is an index pair for $S$.  Therefore, putting $N:=\Cls(N_{2}\setminus N_{0})$ and $N';=\Cls(N_{2}\setminus N_{1})$ and applying Sublemma \ref{Nested_IsoNbd}, we have $A^{\ast}=I(\Cls(N_{2}\setminus N_{1}))$.  
\hfill\textit{QED of Lemma \ref{index_triple}.}

\begin{defn}
 The triple $(N_{2}, N_{1}, N_{0})$ as in Lemma \ref{index_triple} is called an {\em index triple} for the attractor-repeller pair $(A, A^{\ast})$ in $S$.  
\end{defn}

Following Franzosa \cite{Franzosa}, we will generalise this index triple for an attractor-repeller pair to a global notion for a Morse decomposition.  Recall that $\mathcal{M}=\{ m(\pi )\}_{\pi\in P}$ is a Morse decomposition for $(M, F)$.  

\begin{defn}\label{indx_filt}
An {\em index filtration} for $\mathcal{M}$ is a family $\mathcal{N}=\{ N(I)\}_{I\in\mathscr{A}}$ of compact subsets which satisfies the following:
\begin{enumerate}
 \item\label{i_f_1} For any $I\in\mathscr{A}$, $(N(I), N(\emptyset ))$ is an index pair for the attractor $m(I)$.  
 \item  For any $I_{1}, I_{2}\in\mathscr{A}$, $N(I_{1}\cap I_{2})=N(I_{1})\cap N(I_{2})$ and $N(I_{1}\cup I_{2})=N(I_{1})\cup N(I_{2})$.  
\end{enumerate}
\end{defn}
If $\mathcal{N}$ is an index filtration for the Morse decomposition $\mathcal{M}$, for any attractor $m(I)\ (I\in\mathscr{A})$ there exists its index pair in $\mathcal{N}$ (by the condition (\ref{i_f_1}) of Definition \ref{indx_filt}).  Moreover, by (\ref{ip_MS}) of the following proposition, it follows that for any Morse set $m(J)\ (J\in\mathscr{I})$, there exists its index pair in $\mathcal{N}$.  

\begin{prop}\label{w/o_proof}
 Suppose $\mathcal{N}=\{ N(I)\}_{I\in\mathscr{A}}$ is an index filtration for the Morse decomposition $\mathcal{M}$ and $J\in\mathscr{I}$ is given, then we have the following:
\begin{enumerate}
 \item\label{ip_MS} Let $(I, J)$ be a decomposition of an attracting interval $K\in\mathscr{A}$.  Then, $(N(K), N(I))$ is an index pair of the Morse set $m(J)$.  
 \item For $i=1, 2$, suppose $(I_{i}, J)$ is a decomposition of an attracting interval $K_{i}$, then $N(K_{1})\setminus N(I_{1})=N(K_{2})\setminus N(I_{2})$.  
\end{enumerate}
\end{prop}

We omit the proofs of an existence theorem of index filtrations for any Morse decompositions and Proposition \ref{w/o_proof}, since we have a natural, simple index filtration for an NMS-flow (see Section \ref{NMS_RHD}).  For the detail of the existence theorem in the general case, see \cite{Franzosa}.  

\subsection{Homology index braids}\label{homology_index_braid}
Let $S$ be a compact invariant subset and $(A, A^{\ast})$ an attractor-repeller pair in $S$.  Suppose $(N_{2}, N_{1}, N_{0})$ is an index triple of $(A, A^{\ast})$.  Then we have maps of pointed spaces induced from the inclusions
\[
 N_{1}/N_{0}\stackrel{i}{\rightarrow} N_{2}/N_{0}\stackrel{p}{\rightarrow} N_{2}/N_{1} 
\]
and chain maps
\[
 C_{\ast}(N_{1}/N_{0})\stackrel{i_{\ast}}{\rightarrow} C_{\ast}(N_{2}/N_{0})\stackrel{p_{\ast}}{\rightarrow} C_{\ast}(N_{2}/N_{1}).  
\]
On the other hand, we have a short exact sequence 
\[
 0\rightarrow C_{\ast}(N_{1}/N_{0})\stackrel{i_{\ast}}{\rightarrow} C_{\ast}(N_{2}/N_{0})\stackrel{q_{\ast}}{\rightarrow} C_{\ast}(N_{2}/N_{0}, N_{1}/N_{0})\rightarrow 0.  
\]
Since $p_{\ast}\circ i_{\ast}=0$ we have a homomorphism $\rho : C_{\ast}(N_{2}/N_{0}, N_{1}/N_{0})\rightarrow C_{\ast}(N_{2}/N_{1})$ with $\rho\circ q_{\ast}=p_{\ast}$.  
\begin{prop}\label{rel_qtnt}
 For any $m\in\Z_{>0}$ the homomorphism $\rho_{\ast}:H_{m}(N_{2}/N_{0}, N_{1}/N_{0})\rightarrow H_{m}(N_{2}/N_{1})$ induced by $\rho$ is an isomorphism.  
\end{prop}
\begin{proof}
 By identifying $(N_{2}/N_{0})/(N_{1}/N_{0})=N_{2}/N_{1}$, $\rho_{\ast}$ is equal to a homomorphism $\epsilon_{\ast}:H_{\ast}(N_{2}/N_{0}, N_{1}/N_{0})\rightarrow H_{\ast}((N_{2}/N_{0})/(N_{1}/N_{0}))$.  It is shown by Kurland \cite{Kurland} (Proof of Lemma 3.4) that there exists a closed neighbourhood $U$ of $N_{1}/N_0$ in $N_{2}/N_{0}$ such that $N_{1}/N_{0}$ is a weak deformation retract of $U$.  It can be seen that for the pair $(N_{2}/N_{0}, N_{1}/N_{0})$ the natural homomorphism $H_{\ast}(N_{2}/N_{0}, N_{1}/N_{0})\rightarrow H_{\ast}((N_{2}/N_{0})/(N_{1}/N_{0}), (N_{1}/N_{0})/(N_{1}/N_{0}))=H_{\ast}(N_{2}/N_{1}, N_{1}/N_{1})$ is an isomorphism.  Since $H_{\ast}(N_{2}/N_{1}, N_{1}/N_{1})$ is isomorphic to the reduced homology $\widetilde{H}_{\ast}(N_{2}/N_{1})$, the homomorphism $\epsilon_{\ast}(=\rho_{\ast})$ is an isomorphism for $\ast >0$.   
\end{proof}

\begin{remark}
 In our case of NMS-flows, we need not the full generality such as Kurland's result in the above proof (See Section \ref{NMS_RHD}, \ref{def_bdry}).  
\end{remark}

\noindent
Hence from the above short exact sequence we have a long exact sequence of homology 
\[
 \rightarrow H_{m+1}(N_{1}/N_{0})\rightarrow H_{m+1}(N_{2}/N_{0})\stackrel{\rho_{\ast}}{\rightarrow}H_{m+1}(N_{2}/N_{1})\stackrel{\partial}{\rightarrow}H_{m}(N_{1}/N_{0})\rightarrow
\]
where we replace $H_{m+1}((N_{2}/N_{0})/(N_{1}/N_{0})$ with $H_{m+1}(N_{2}/N_{1})$ by Proposition \ref{rel_qtnt}.  As homology indices we have a long exact sequence 
\begin{equation*}
 \vcenter{
  \xymatrix@C=30pt{
   \mathscr{H}_{\ast}(A) \ar[rr]^{i(A, S)} & & \mathscr{H}_{\ast}(S) \ar[dl]^{p(S, A^{\ast})} \\
     & \mathscr{H}_{\ast}(A^{\ast}) \ar[lu]^{\partial (A^{\ast}, A)} & 
}
}
\end{equation*}
where all arrows are homomorphism of graded modules and $\mathrm{deg}(i(A, S))=\mathrm{deg}(p(S, A^{\ast}))=0$ and $\mathrm{deg}(\partial (A^{\ast}, A))=-1$.  
We call this long exact sequence the {\em homology index sequence} of the attractor-repeller pair $(A, A^{\ast})$.  

\begin{remark}
 Up to isomorphisms, we can write $C_{\ast}(A)=C_{\ast}(N_{1}/N_{0}), C_{\ast}(S)=C_{\ast}(N_{2}/N_{0}), C_{\ast}(A^{\ast})=C_{\ast}(N_{2}/N_{1})$ and $H_{\ast}(A)=H_{\ast}(N_{1}/N_{0}), H_{\ast}(S)=H_{\ast}(N_{2}/N_{0}), H_{\ast}(A^{\ast})=H_{\ast}(N_{2}/N_{1})$.  
\end{remark}

In the following in this section, we fix a Morse decomposition $\mathcal{M}=\{ m(\pi)\}_{\pi\in P}$ for $(M, F)$ and an index filtration $\mathcal{N}=\{ N(I)\}_{I\in\mathscr{A}}$ for $\mathcal{M}$.  
Suppose an arbitrary $(I, J)\in\mathscr{I}_{2}$ is given.  Define $H$ to be the set of all $\pi\in P\setminus IJ$ such that there exists $\pi'\in IJ$ with $\pi <\pi'$.  Then $HI$ and $HIJ$ are attracting intervals and it can be seen that the index filtration $\mathcal{N}$ provides an index triple $(N(HIJ), N(HI), N(H))$ of $(m(I), m(J))$.  Note that in $m(IJ)$ the complementary repeller $m(I)^{\ast}$ to $m(I)$ is equal to $m(J)$.  Now, as in the above argument we have a sequence of chain maps
\[
 C_{\ast}^{\mathcal{N}}(m(I))\rightarrow C_{\ast}^{\mathcal{N}}(m(IJ))\rightarrow C_{\ast}^{\mathcal{N}}(m(J))
\]
and the homology index sequence of $(m(I), m(J))$
\begin{equation*}
 \vcenter{
  \xymatrix@C=30pt{
   \mathscr{H}_{\ast}^{\mathcal{N}}(m(I)) \ar[rr]^{i(m(I), m(IJ))} & & \mathscr{H}_{\ast}^{\mathcal{N}}(m(IJ)) \ar[dl]^{p(m(IJ), m(J))} \\
     & \mathscr{H}_{\ast}^{\mathcal{N}}(m(J)) \ar[lu]^{\partial (m(J), m(I))} & 
}
}
\end{equation*}
where we add a superscript $\mathcal{N}$ in order to note the role of fixed index filtration.  Also, we may call it the homology index sequence of $(I, J)$ instead of $(m(I), m(J))$.  

Next, let $(I, J, K)\in\mathscr{I}_{3}$ be an arbitrary adjacent triple.  Define $H$ to be the set of all $\pi\in P\setminus IJK$ such that there exists $\pi'\in IJK$ with $\pi <\pi'$.  Then it can be easily seen that $H\in\mathscr{A}$ and that $(H, I, J, K)\in\mathscr{I}_{4}$.  Hence we have $HI, HIJ, HIJK\in\mathscr{A}$.  Set $N_{0}=N(H), N_{1}=N(HI), N_{2}=N(HIJ)$ and $N_{3}=N(HIJK)$, for the simplicity.  Then the following diagram of index spaces and maps induced by inclusions is commutative:

\begin{equation*}
\vcenter{
\xymatrix@C=30pt{
N_{1}/N_{0} \ar[rd]^{i} \ar@/_50pt/[dd]_{i} &  & \\
  & N_{2}/N_{0} \ar[rd]^{p} \ar[ld]_{i} & \\
N_{3}/N_{0} \ar[rd]^{p} \ar@/_50pt/[dd]_{p} & & N_{2}/N_{1} \ar[ld]_{i} \\
  & N_{3}/N_{1} \ar[ld]_{p} & \\
N_{3}/N_{2} & & 
}
}
\end{equation*}

\noindent
Taking the chain complexes, we have the following commutative diagram: 

\begin{equation*}
\vcenter{
\xymatrix@C=50pt{
C^{\mathcal{N}}_{\ast}(m(I)) \ar[rd]^{i} \ar@/_60pt/[dd]_{i} &  & \\
  & C^{\mathcal{N}}_{\ast}(m(IJ)) \ar[rd]^{p} \ar[ld]_{i} & \\
C^{\mathcal{N}}_{\ast}(m(IJK)) \ar[rd]^{p} \ar@/_50pt/[dd]_{p} & & C^{\mathcal{N}}_{\ast}(m(J)) \ar[ld]_{i} \\
  & C^{\mathcal{N}}_{\ast}(m(JK)) \ar[ld]_{p} & \\
C^{\mathcal{N}}_{\ast}(m(K)) & & 
}
}
\end{equation*}

\noindent
This diagram is a combination of four short exact sequences of chain complexes.  Moreover, applying the homology functor, we have the following diagram of the homology indices:

\begin{equation*} 
\vcenter{
\xymatrix@C=50pt{
 \ar@{-->}@/_30pt/[d] & \mbox{\hspace{40pt}}\ar@{-->}[ld] \ar@{-->}[rd]  & \ar@{-->}@/^30pt/[d] \\
\mathscr{H}^{\mathcal{N}}_{\ast}(m(I)) \ar[rd]^{i} \ar@/_50pt/[dd]_{i} & & \mathscr{H}^{\mathcal{N}}_{\ast}(m(K)) \ar[ld]_{\partial} \ar@/^50pt/[dd]^{\partial} \\
 & \mathscr{H}^{\mathcal{N}}_{\ast}(m(IJ)) \ar[ld]_{i} \ar[rd]^{p}  & \\
\mathscr{H}^{\mathcal{N}}_{\ast}(m(IJK)) \ar@/_50pt/[dd]_{p} \ar[rd]^{p} & & \mathscr{H}^{\mathcal{N}}_{\ast}(m(J)) \ar[ld]_{i} \ar@/^50pt/[dd]^{\partial}\\
 & \mathscr{H}^{\mathcal{N}}_{\ast}(m(JK)) \ar[ld]_{p} \ar[rd]^{\partial} & \\
\mathscr{H}^{\mathcal{N}}_{\ast}(m(K)) \ar@/_50pt/[dd]_{\partial} \ar[rd]^{\partial} & & \mathscr{H}^{\mathcal{N}}_{\ast}(m(I)) \ar[ld]_{i} \ar@/^50pt/[dd]^{i}  \\
 & \mathscr{H}^{\mathcal{N}}_{\ast}(m(IJ)) \ar[ld]_{p} \ar[rd]^{i}& \\
\mathscr{H}^{\mathcal{N}}_{\ast}(m(J)) \ar@{-->}@/_30pt/[d] \ar@{-->}[rd] & & \mathscr{H}^{\mathcal{N}}_{\ast}(m(IJK)) \ar@{-->}[ld] \ar@{-->}@/^30pt/[d] \\
 & \mbox{\hspace{70pt}}  & 
}
} 
\end{equation*}
\noindent
This is a combination of the homology index sequences of $(m(I), m(J)), (m(I), m(JK)), (m(J), m(K))$ and $(m(IJ), m(K))$.  We call this diagram {\em homology index braid} of $(I, J, K)$.

\section{Nonsingular Morse-Smale flows and round handle decompositions}\label{NMS_RHD} 
Let $M$ be a compact $n$-manifold, possibly with boundary.  A 
flow $F$ on $M$ is called a {\em Morse-Smale flow} if it is generated by a vector field $X$ which satisfies the following conditions:
\begin{enumerate}
\item 
$X$ has a finite number of critical elements (singular points and periodic orbits), all of which are hyperbolic,
\item
for any critical elements $\sigma_{1}, \sigma_{2}$ of $X$, the unstable submanifold $W^{u}(\sigma_{1})$ is transversal to the stable submanifold $W^{s}(\sigma_{2})$,
\item
the nonwandering set $\Omega (X)$ of $X$ is equal to the union of the critical elements of $X$.  
\end{enumerate}
For the detailed property of Morse-Smale vector fields and flows, we refer the readers to \cite{Palis-Melo}.  
Recall that we concern a one-dimensional oriented foliation, an {\em NMS-foliation}, which consists of all the orbits of an NMS-flow as leaves.  

Suppose that a union $\partial_{-}M$ of connected components of the boundary $\partial M$ is specified.  We call such a pair $(M, \partial_{-}M)$ a {\em manifold pair} in this paper.  
A vector field on a manifold pair $(M, \partial_{-}M)$ is assumed to be transversal to $\partial M$ and pointing outward on $\partial_{-}M$ and inward on $\partial_{+}M := \partial M\setminus\partial_{-}M$, unless otherwise stated.  Note that the necessary and sufficient condition for the existence of a {\em nonsingular} vector field on $(M, \partial_{-}M)$ is vanishing of its relative Euler characteristic $\chi (M, \partial_{-}M)$. 

Next, we define a round handle decomposition of $(M, \partial_{-}M)$.  Let $E^{k}$ denote a $k$-disk bundle over the circle $S^{1}$ where $k$ is an integer with $0\leq k\leq n-1$.  Set $R^{k} := E^{k}\oplus E^{n-k-1}$, $\partial_{-}R^{k} := \partial E^{k}\times_{S^{1}}E^{n-k-1}$ and $\partial_{+}R^{k} := E^{k}\times_{S^{1}}\partial E^{n-k-1}$.  We call $R^{k}$ or $(R^{k}, \partial_{-}R^{k})$ an {\em $n$-dimensional round handle of index $k$} or {\em $n$-dimensional round $k$-handle}. The zero-section of a round $k$-handle is called the {\em core circle} of the round handle.

Note that $R^{k}$ is orientable if and only if the both of the disk bundles $E^{k}$ and $E^{n-k-1}$ are trivial or non-trivial simultaneously.  We say $R^{k}$ is a {\em trivial} round handle if the both of $E^{k}$ and $E^{n-k-1}$ are trivial and {\em twisted} otherwise.  Also, if at least $E^{k}$ is trivial, then $R^{k}$ is said to be {\em quasi-trivial}.  

Let $(W, \partial_{-}W)$ be an $n$-dimensional manifold pair.  
For an embedding $\alpha :\partial_{-}R^{k}\rightarrow \partial_{+}W$, we define the identification space $W\cup_{\alpha}R^{k}$, the quotient space of the disjoint union $W\sqcup R^{k}$ identified by $\alpha$.  The resulting manifold is also denoted by $W +_{\alpha}R^{k}$ or simply by $W + R^{k}$.  Set $\partial_{-}(W +_{\alpha} R^{k}) = \partial_{-}W$ and note that $\partial_{+}(W +_{\alpha} R^{k}) = (\partial_{+}W\setminus\mathrm{Im}\alpha)\cup\partial_{+}R^{k}$.  We say that $W +_{\alpha}R^{k}$ is obtained from $W$ by {\em attaching a round $k$-handle $R^{k}$}.   
A {\em round handle decomposition} (RHD for short) of a manifold pair $(M, \partial_{-}M)$ is a representation 
\[
 M = (\partial_{-}M\times [0, 1])+R^{k_{1}}+\cdots +R^{k_{m}}
\]
which denotes successive operations of attaching round handles.  We say an RHD is {\em simple} if all the round handles appearing in the RHD are trivial.  
\vspace{2ex}

The following propositions are foundations of our research. 

\begin{prop}[Asimov \cite{Asimov_RHD}]
Let $(M, \partial_{-}M)$ is a manifold pair with $\chi (M, \partial_{-}M) = 0$.  Then, $(M, \partial_{-}M)$ has an RHD if $\mathrm{dim}M \neq 3$.  Moreover, if $M$ is not a M\"{obius band} $($and $\mathrm{dim}\neq 3)$, then $(M, \partial_{-}M)$ has a simple RHD.
\end{prop}

\begin{prop}[Morgan \cite{Morgan}]
 Suppose $(M, \partial_{-}M)$ is a manifold pair and $\mathrm{dim}M = 3$.  Then, $(M, \partial_{-}M)$ has an RHD if and only if $M$ is a graph manifold.  
\end{prop}

\begin{prop}[Asimov \cite{Asimov_RHD}, \cite{Asimov_StrStable}]
If a  manifold pair $(M, \partial_{-}M)$ has an RHD, then $(M, \partial_{-}M)$ has an NMS-flow whose periodic orbits are exactly the core circles of the RHD.  Moreover, any nonsingular vector field on $(M, \partial_{-}M)$ is homotopic to an NMS vector field through nonsingular vector fields.  
\end{prop}

\begin{prop}[Morgan \cite{Morgan}]\label{NMS2RHD}
If a manifold pair $(M, \partial_{-}M)$ has an NMS-flow, then $(M, \partial_{-}M)$ has an RHD whose core circles are exactly the periodic orbits of the flow.
\end{prop}
\vspace{2ex}

\noindent
Here, the correspondence between an NMS-flow and an RHD is obtained as follows.  A round handle $(R^{k}, \partial_{-}R^{k})$ may be considered as a manifold pair with a flow.  
A periodic orbit of index $k$ is just the core circle of a round $k$-handle $R^{k}=E^{k}\oplus E^{n-k-1}$ and its unstable (resp. stable) submanifold intersects $R^{k}$ at $E^{k}\oplus 0$ (resp. $0\oplus E^{n-k-1}$).  This round handle $R^{k}$ is said to be {\em associated} with the periodic orbit.  It is crucial that an associated round handle $(R^{k}, \partial_{-}R^{k})$ is an index pair of the periodic orbit.  

An RHD can be characterised by {\em round Morse functions} or {\em Morse-Bott functions} with one-dimensional critical submanifolds.  Such a viewpoint shows faintly that an NMS-flow has gradient-like property (cf. \cite{mytop}, \cite{Banyaga-Hurtubise}, \cite{Asimov_RHD}).  

Suppose that an NMS-flow $F$ is given on a manifold pair $(M, \partial_{-}M)$ and let $c_{1}$ and $c_{2}$ be periodic orbits of $F$.  Suppose the index of $c_{i}$ is $k_{i}$ ($i = 1, 2$).  Since the unstable submanifold $W^{u}(c_{1})$ of $c_{1}$ and the stable submanifold $W^{s}(c_{2})$ of $c_{2}$ are transversal to each other, we have $W^{u}(c_{1})\cap W^{s}(c_{2}) = \emptyset$ if $k_{1} < k_{2}$.  Note that (un)stable submanifolds are invariant.

\section{Boundary homomorphisms for general NMS-flows}\label{def_bdry}
In this section, we define a homology of an arbitrary NMS-foliation.  
Let $M$ be a closed $n$-manifold and $F$ an NMS-flow.  The NMS-foliation defined by $F$ is denoted by $\mathcal{F}$.  Let $\mathcal{M}=\{ m(\pi)\}_{\pi\in P}$ be the set of all periodic orbits of $F$.  Here the index set $P$ is equal to $\mathcal{M}$ as sets, and it has a natural partial order $\leq_{F}$ defined as follows: for $\pi, \pi'\in P$, $\pi <_{F}\pi'$ if there is a sequence $\pi =\pi_{0}, \pi_{1}, \ldots, \pi_{m}=\pi'\in P$ such that $\pi_{i}\ne\pi_{j}$ if $i\ne j$ and $C(m(\pi_{j}), m(\pi_{j-1}))\ne\emptyset\ (j=1, 2, \ldots, m)$.  It is easily seen that $\leq_{F}$ is a partial order on $P$.  We call this order the {\em flow-order} with respect to $F$.  In the rest of the paper, the order on $P$ is the flow-order with respect to $F$ and we simply write $\leq$ instead of $\leq_{F}$.  
Also, let $\mathcal{R}=\{ R(\pi)\ |\ \pi\in P\}$ be the RHD associated with $F$, where $R(\pi)$ denotes the round handle associated with $\pi$.  

The following lemmata can be easily proven:

\begin{lem}
 $\mathcal{M}$ is a Morse decomposition for $(M, F)$.  
\end{lem}

\begin{lem}
 $\mathcal{R}$ is an index filtration for $\mathcal{M}$.  In particular, for any $\pi\in P$, $(R(\pi), \partial_{-}R(\pi))$ is an index pair for $m(\pi)$.  
\end{lem}

\begin{lem}
 For any $\pi\in P$ of index $k$, we have
\[
 H_{m}(R(\pi), \partial_{-}R(\pi);\ZZ)\cong
\left\{
\begin{array}{cl}
 \ZZ & \text{if } m=k \text{ or } k+1 \\
 0  & \text{otherwise},  
\end{array}
\right.  
\]
\noindent
and
if $R(\pi)$ is a quasi-trivial round handle, then
\[
 H_{m}(R(\pi), \partial_{-}R(\pi);\Z)\cong
\left\{
\begin{array}{cl}
 \Z & \text{if }m=k \text{ or }k+1 \\
 0  & \text{otherwise}.  
\end{array}
\right.
\]
\end{lem}
\noindent
Since $(R(\pi), \partial_{-}(\pi))$ is an index pair of $m(\pi)$, the homology index of $\pi$ of index $k$ is isomorphic to $H_{k}(R(\pi), \partial_{-}R(\pi))\oplus H_{k+1}(R(\pi), \partial_{-}R(\pi))$.

From now on, we identify $\mathcal{M}$ and $P$ as sets: $\mathcal{M}=P, m(\pi)=\pi$.  For $k\in\Z,\ 0\leq k\leq n-1$, set $P^{k}=\{\pi\in P\ |\ \mathrm{ind}(\pi)=k\}$, where $\mathrm{ind}(\pi)$ denotes the index of the hyperbolic periodic orbit $\pi$.  We consider a periodic orbit of index $k$ as a generator of our $k$-chain groups.  For any $\pi^{k}\in P^{k}$, set
\[
\begin{array}{lcl}
 C^{\mathcal{F}}_{\langle k+1\rangle}(\pi^{k}) & := & C^{\mathcal{R}}_{\langle k+1\rangle}(\pi^{k})=H_{k+1}(R^{k}, \partial_{-}R^{k}) \\
 C^{\mathcal{F}}_{\langle k\rangle}(\pi^{k}) & := & C^{\mathcal{R}}_{\langle k\rangle}(\pi^{k})=H_{k}(R^{k}, \partial_{-}R^{k}) \\
 C_{k}(\mathcal{F}) & := & \left(\oplus_{\pi\in P^{k}}C^{\mathcal{F}}_{\langle k+1\rangle}(\pi)\right)\oplus\left(\oplus_{\pi\in P^{k}}C^{\mathcal{F}}_{\langle k\rangle}(\pi)\right). 
\end{array}
\]
Note that $\langle\ell\rangle$ denotes the dimension of homology of the round handle not its index.  
Here and in the following, the coefficient of homology groups depends on the type of the associated RHD.  That is, it is $\Z$ if its RHD is simple, and $\ZZ$ if its RHD contains a twisted round handle.  The elements of $C^{\mathcal{F}}_{\langle k+1\rangle}(\pi^{k})$ and $C^{\mathcal{F}}_{\langle k\rangle}(\pi^{k})$ are called {\em fundamental $k$-chains} and {\em sectional $k$-chains} of $\mathcal{F}$ respectively.  

For any $\pi\in P$, we define an interval $\mathrm{Blw}(\pi)$ as the set of all elements each of which is directly below $\pi$.  That is, we set
\[
\begin{array}{lcl}
 \mathrm{Blw}(\pi) & := & \{\pi'\in P\ |\ \pi'<\pi \text{ and there is no }\pi''\in P \text{ such that }\pi'<\pi''<\pi\} \\
 \mathrm{Blw}(\pi)^{k} & := & \mathrm{Blw}(\pi)\cap P^{k}.  
\end{array}
\]
We will define a boundary homomorphism $\partial^{\mathcal{F}}:\mathscr{H}^{\mathcal{R}}_{\ast}(\{\pi\})\rightarrow\mathscr{H}^{\mathcal{R}}_{\ast}(\mathrm{Blw}(\pi))$ and show that $\partial^{\mathcal{F}}\circ\partial^{\mathcal{F}}=0$. 

We fix an arbitrary $\xi\in P^{k+1}$.  Set $K:=\{\xi\}$ and $J:=\mathrm{Blw}(\xi)$.  It is clear that $J, K\in\mathscr{I}, J\cap K=\emptyset$ and any two elements of $J$ are non-comparable.  Next, set $I:=\mathrm{Itv}(\cup_{\eta\in J}\mathrm{Blw}(\eta))$, where $\mathrm{Itv}(A)$ of $A\subset P$ denotes the {\em intervalisation} of $A$, i.e., the smallest interval which contains $A$.  Moreover, let $H$ be the set of all $\pi\in P\setminus IJK$ such that there exists $\pi'\in IJK$ with $\pi <\pi'$.  Then it is clear that $H, I, J, K\in\mathscr{I}$ and $I\cap JK=\emptyset$.  
\begin{note}
 In order only to define $\partial^{\mathcal{F}}$, we need not to set the interval $I$ and it is enough to set $H':=\{\pi\in P\setminus JK\ |\ \text{there is }\pi'\in JK\text{ such that }\pi<\pi'\}$ instead.  It is easily seen that $H'=HI$ and the following argument to define $\partial^{\mathcal{F}}$ holds if we replace $HI$ with $H'$.  
\end{note}
We will show the properties of $H, I, J$ and $K$ in the following.  

\setcounter{claim}{0}
\begin{claim}
 We have $HIJK=\{\pi\in P\ |\ \pi\leq\xi\}$.  In particular, $HIJK\in\mathscr{A}$.  
\end{claim}
\begin{proof}
 Suppose $\pi\in HIJK$ is given.  If $\pi\in JK$, then clearly $\pi\leq\xi$.  Suppose $\pi\in I$.  We have to show in the case $\pi\not\in\cup_{\eta\in J}\mathrm{Blw}(\eta)$.  By the minimality of intervalisation, there exist $\pi', \pi''\in\cup_{\eta\in J}\mathrm{Blw}(\eta)$ such that $\pi'<\pi<\pi''$.  Then since $\pi''<\xi$ we have $\pi<\xi$.  Lastly, assume $\pi\in H$. Then there exists $\pi'\in IJK$ such that $\pi<\pi'$.  As we showed, $\pi'\leq\xi$, thus $\pi<\xi$.  

On the contrary, suppose $\pi\in P$ with $\pi\leq\xi$ is given.  We may assume $\pi<\xi$.  If $\pi\not\in IJK$, then by the definition of $H$, we have $\pi\in H$.    
\end{proof}
Next, we show that $(H, I, J, K)$ satisfies the order condition of adjacent 4-tuple, that is, for disjoint two intervals $A$ and $B$ appearing in this order in a tuple, if $a\in A$ and $b\in B$, then $b\nless a$.  We write $A\prec B$ if this condition is satisfied.  First, $H\prec K, I\prec K, J\prec K$ are clear.  
\begin{claim}
 $H\prec J, I\prec J$.  
\end{claim}

\begin{proof}
 Assume that there exist $\pi\in HI$ and $\pi_{J}\in J$ such that $\pi_{J}<\pi$.  Then, there exists $\eta\in IJK$ such that $\pi <\eta$.  Precisely, if $\pi\in H$ then $\eta\in IJK$ and if $\pi\in I$ then $\eta\in J$.  
Thus we have $\pi_{J}<\pi <\eta\leq\xi$, and this contradicts to $\pi_{J}\in J$.  
\end{proof}
\begin{claim}
 $H\prec I$.  
\end{claim}
\begin{proof}
 Assume there exist $\pi_{I}\in I$ and $\pi_{H}\in H$ such that $\pi_{I}<\pi_{H}$.  Then, we claim there exists $\pi\in I$ such that $\pi_{H}<\pi$.  First, there exists $\eta\in IJK$ such that $\pi_{H}<\eta$ by the definition of $H$.  If $\eta\not\in I$, there exists $\pi\in I$ such that $\pi_{H} <\pi <\eta$ since $\leq$ is the flow-order.  Then, we have $\pi_{I}<\pi_{H}<\pi$ and $\pi_{I}, \pi\in I$, which implies $\pi_{H}\in I$, a contradiction to $H\cap I=\emptyset$.   
\end{proof}
Thus, we have $(H, I, J, K)\in\mathscr{I}_{4}$. Also, we have $H\in\mathscr{A}$, therefore it is clear that $HI, HIJ, HIJK\in\mathscr{A}$.  
We will define $\partial^{\mathcal{F}}$ by applying the homology index sequence of $(J, K)$ and show that $\partial^{\mathcal{F}}\circ\partial^{\mathcal{F}}=0$ by applying the homology index braid of $(I, J, K)$.  
\vspace{2ex}

\noindent
\textbf{Definition of }$\boldsymbol{\partial^{\mathcal{F}}}$

Now, we have maps induced from inclusions between index spaces:
\[
 R(HIJ)/R(HI)\rightarrow R(HIJK)/R(HI)\rightarrow R(HIJK)/R(HIJ)
\]
which induces the sequence of chain complexes:
\[
 C^{\mathcal{R}}_{\ast}(J)\rightarrow C^{\mathcal{R}}_{\ast}(JK)\rightarrow C^{\mathcal{R}}_{\ast}(K).  
\]
Applying the homology functor, we have the homology index sequence as in the argument in Section \ref{MD_CI_HIB}:
\begin{equation*}
 \vcenter{
  \xymatrix@C=30pt{
   \mathscr
{H}^{\mathcal{R}}_{\ast}(J) \ar[rr]^{i(J, JK)} & & \mathscr{H}^{\mathcal{R}}_{\ast}(JK) \ar[dl]^{p(JK, K)} \\
     & \mathscr{H}^{\mathcal{R}}_{\ast}(K) \ar[lu]^{\partial (K, J)} & 
}
}
\end{equation*}
where $C^{\mathcal{R}}_{\ast}(m(J))$ and $\mathscr{H}^{\mathcal{R}}_{\ast}(m(I))$ are abbreviated to $C^{\mathcal{R}}_{\ast}(J)$ and $\mathscr{H}^{\mathcal{R}}_{\ast}(I)$ respectively, etc.  Since $(R(\xi), \partial_{-}R(\xi))$ is an index pair of $\xi$, we have
\begin{eqnarray*}
 \mathscr{H}^{\mathcal{R}}_{\ast}(K) & = & H_{\ast}(R(HIJK)/R(HIJ)) \\
  & \cong & H_{\ast}(R(K), \partial_{-}R(K))\hspace{2em} (\ast>0) \\
  & = & H_{\ast}(R(\xi), \partial_{-}R(\xi)) \\
  & = & C^{\mathcal{F}}_{\langle k+2\rangle}(\xi )\oplus C^{\mathcal{F}}_{\langle k+1\rangle}(\xi )
\end{eqnarray*}
and also we have
\begin{eqnarray*}
 \mathscr{H}^{\mathcal{R}}_{\ast}(J) & = & H_{\ast}(R(HIJ)/R(HI)) \\
  & \cong & H_{\ast}(R(J)/\partial_{-}R(J)) \\
  & \cong & H_{\ast}\left(\sqcup_{\pi\in J}R(\pi)/\sqcup_{\pi\in J}\partial_{-}R(\pi)\right) \\
  & \cong & \oplus_{\pi\in J}H_{\ast}(R(\pi), \partial_{-}R(\pi)) \hspace{2em} (\ast>0) \\
  & = & \oplus_{\pi\in J}\left( C^{\mathcal{F}}_{\langle \ell +1\rangle}(\pi)\oplus C^{\mathcal{F}}_{\langle \ell\rangle}(\pi)\right)\hspace{1.5em}(\ell =\mathrm{ind}(\pi)\leq k) 
\end{eqnarray*}
Let $[\xi]_{f}$ and $[\xi]_{s}$ be bases of $C^{\mathcal{F}}_{\langle k+2\rangle}(\xi)$ and $C^{\mathcal{F}}_{\langle k+1\rangle}(\xi)$ respectively, and we define
\begin{eqnarray*}
 \partial^{\mathcal{F}}([\xi]_{f}) & = & \partial (K, J)([\xi]_{f})\hspace{1.5em}(\in\mathscr{H}^{\mathcal{R}}_{k+1}(J)) \\ 
 \partial^{\mathcal{F}}([\xi]_{s}) & = & \partial (K, J)([\xi]_{s})\hspace{1.5em}(\in\mathscr{H}^{\mathcal{R}}_{k}(J)).  
\end{eqnarray*}
Thus, we have a boundary homomorphism
\[
 \partial^{\mathcal{F}}:C_{k+1}(\mathcal{F})\rightarrow\oplus_{\ell\leq k+1}C_{\ell}(\mathcal{F}).  
\]
\begin{remark}
 As for the dimension of the target chains (= the index of the target periodic orbits), there are some cases.  We will specify the cases later.  
\end{remark}

\noindent
\textbf{Assertion}
 $\partial^{\mathcal{F}}\circ\partial^{\mathcal{F}}=0$.  

To prove the assertion, we consider the homology index braid determined from $(H, I, J, K)\in\mathscr{I}_{4}$ and $H, HI, HIJ, HIJK\in\mathscr{A}$.  As we discussed it in the Subsection \ref{homology_index_braid}, we have the following commutative diagrams:
\begin{equation*}
\vcenter{
\xymatrix@C=30pt{
R(HI)/R(H) \ar[rd]^{i} \ar@/_50pt/[dd]_{i} &  & \\
  & R(HIJ)/R(H) \ar[rd]^{p} \ar[ld]_{i} & \\
R(HIJK)/R(H) \ar[rd]^{p} \ar@/_50pt/[dd]_{p} & & R(HIJ)/R(HI) \ar[ld]_{i} \\
  & R(HIJK)/R(HI) \ar[ld]_{p} & \\
R(HIJK)/R(HIJ) & & 
}
}
\end{equation*}
\begin{equation*}
\vcenter{
\xymatrix@C=50pt{
C^{\mathcal{R}}_{\ast}(I) \ar[rd]^{i} \ar@/_60pt/[dd]_{i} &  & \\
  & C^{\mathcal{R}}_{\ast}(IJ) \ar[rd]^{p} \ar[ld]_{i} & \\
C^{\mathcal{R}}_{\ast}(IJK) \ar[rd]^{p} \ar@/_50pt/[dd]_{p} & & C^{\mathcal{R}}_{\ast}(J) \ar[ld]_{i} \\
  & C^{\mathcal{R}}_{\ast}(JK) \ar[ld]_{p} & \\
C^{\mathcal{R}}_{\ast}(K) & & 
}
}
\end{equation*}
This diagram induces the following homology index braid:
\begin{equation*} 
\vcenter{
\xymatrix@C=50pt{
 \mbox{\hspace{70pt}} \ar@{-->}@/_30pt/[dd] \ar@{-->}[rd] & & \mbox{\hspace{70pt}} \ar@{-->}[ld]   \ar@{-->}@/^30pt/[dd] \\
  & \mathscr{H}^{\mathcal{R}}_{\ast}(JK) \ar[ld]_{\partial} \ar[rd]^{p} \\
\mathscr{H}^{\mathcal{R}}_{\ast}(I) \ar[rd]^{i} \ar@/_50pt/[dd]_{i} & & \mathscr{H}^{\mathcal{R}}_{\ast}(K) \ar[ld]_{\partial} \ar@/^50pt/[dd]^{\partial} \\
 & \mathscr{H}^{\mathcal{R}}_{\ast}(IJ) \ar[ld]_{i} \ar[rd]^{p}  & \\
\mathscr{H}^{\mathcal{R}}_{\ast}(IJK) \ar@/_50pt/[dd]_{p} \ar[rd]^{p} & & \mathscr{H}^{\mathcal{R}}_{\ast}(J) \ar[ld]_{i} \ar@/^50pt/[dd]^{\partial}\\
 & \mathscr{H}^{\mathcal{R}}_{\ast}(JK) \ar[ld]_{p} \ar[rd]^{\partial} & \\
\mathscr{H}^{\mathcal{R}}_{\ast}(K) \ar@/_50pt/[dd]_{\partial} \ar[rd]^{\partial} & & \mathscr{H}^{\mathcal{R}}_{\ast}(I) \ar[ld]_{i} \ar@/^50pt/[dd]^{i}  \\
 & \mathscr{H}^{\mathcal{R}}_{\ast}(IJ) \ar[ld]_{p} \ar[rd]^{i}& \\
\mathscr{H}^{\mathcal{R}}_{\ast}(J) \ar@{-->}@/_30pt/[d] \ar@{-->}[rd] & & \mathscr{H}^{\mathcal{R}}_{\ast}(IJK) \ar@{-->}[ld] \ar@{-->}@/^30pt/[d] \\
 & \mbox{\hspace{70pt}}  & 
}
} 
\end{equation*}
We write $[\xi]_{c}\in C^{\mathcal{F}}_{\langle\ell\rangle}(\xi)=\mathscr{H}^{\mathcal{R}}_{\ell}(K)$, where $c=f$ if $\ell =k+2$ and $c=s$ if $\ell =k+1$.  Then by the commutativity of this homology index braid, we have
\begin{eqnarray*}
 \partial^{\mathcal{F}}(\partial^{\mathcal{F}}([\xi]_{c})) & = & \partial^{\mathcal{F}}(\partial (K, J)([\xi]_{c})) \\
 & = & \partial (J, I)\circ\partial (K, J)([\xi]_{c}) \\
 & = & (\partial (JK, I)\circ i(J, JK))\circ\partial (K, J)([\xi]_{c}) \\
 & = & 0
\end{eqnarray*}
This completes the proof of Assertion.  
\vspace{2ex}

Then we define the homology group $H_{\ast}(\mathcal{F})$ as follows.  Recall that $C_{k}(\mathcal{F}) = \left(\oplus_{\pi\in P^{k}}C^{\mathcal{F}}_{\langle k+1\rangle}(\pi)\right)\oplus\left(\oplus_{\pi\in P^{k}}C^{\mathcal{F}}_{\langle k\rangle}(\pi)\right)$.  We set
\begin{eqnarray*}
 C^{f}_{k}(\mathcal{F}) & := & \oplus_{\pi\in P^{k}}C_{\langle k+1\rangle}(\mathcal{F}) \\
 C^{s}_{k}(\mathcal{F}) & := & \oplus_{\pi\in P^{k}}C_{\langle k\rangle}(\mathcal{F}).  
\end{eqnarray*}
Suppose $\pi\in P^{k}$ and $\eta\in\mathrm{Blw}(\pi)$ are given. Then the dimension (= index) of a non-trivial component at the target $\eta$ of the boundary homomorphism $\partial^{\mathcal{F}}(\pi)$ depends on the types (fundamental or sectional) of the source chains of $\pi$, i.e., the homology dimension h-dim($\pi$)(=$\langle k\rangle$ or $\langle k+1\rangle$) as follows.
\[
\begin{array}{l|l}
 \text{h-dim}(\pi) & \mathrm{ind}(\eta) \\ \hline
 \langle k+1\rangle\ \mathrm{(fund.)} & k-1\ \mathrm{(fund.)} \\
 \langle k+1\rangle\ \mathrm{(fund.)} & k\ \mathrm{(sect.)} \\
 \langle k\rangle\ \mathrm{(sect.)} & k-1\ \mathrm{(sect.)} \\
 \langle k\rangle\ \mathrm{(sect.)} & k-2\ \mathrm{(fund.)} 
\end{array} 
\]
Thus we may write
\begin{eqnarray*}
\partial^{f}_{k} & : & C^{f}_{k}(\mathcal{F})\rightarrow C^{f}_{k-1}(\mathcal{F})\oplus C^{s}_{k}(\mathcal{F}) \\
 \partial^{s}_{k} & : & C^{s}_{k}(\mathcal{F})\rightarrow C^{f}_{k-2}(\mathcal{F})\oplus C^{s}_{k-1}(\mathcal{F}) \\
\partial^{\mathcal{F}}_{k} & = & \partial^{f}_{k}\oplus\partial^{s}_{k}, 
\end{eqnarray*}
\noindent
and hence the boundaries and cycles are denoted by
\begin{eqnarray*}
 B^{f}_{k}(\mathcal{F}) & := & (\mathrm{Im}(\partial^{f}_{k+1})\oplus\mathrm{Im}(\partial^{s}_{k+2}))\cap C^{f}_{k}(\mathcal{F}) \\ 
 B^{s}_{k}(\mathcal{F}) & := & (\mathrm{Im}(\partial^{f}_{k})\oplus\mathrm{Im}(\partial^{s}_{k+1}))\cap C^{s}_{k}(\mathcal{F}) \\
 B_{k}(\mathcal{F}) & := & B^{f}_{k}(\mathcal{F})\oplus B^{s}_{k}(\mathcal{F}) \\
 Z^{f}_{k}({\mathcal{F}}) & := & \mathrm{Ker}(\partial^{f}_{k}) \\
 Z^{s}_{k}(\mathcal{F}) & := & \mathrm{Ker}(\partial^{s}_{k}) \\
 Z_{k}(\mathcal{F}) & := & Z^{f}_{k}(\mathcal{F})\oplus Z^{s}_{k}(\mathcal{F})
\end{eqnarray*}
and define
\[
 H_{k}(\mathcal{F}) := Z_{k}(\mathcal{F})/B_{k}(\mathcal{F}).  
\]
Clearly, we have
\begin{eqnarray*}
 H_{k}(\mathcal{F}) & = & H^{f}_{k}(\mathcal{F})\oplus H^{s}_{k}(\mathcal{F}) \\
 & = & (Z^{f}_{k}(\mathcal{F})/B^{f}_{k}(\mathcal{F}))\oplus (Z^{s}_{k}(\mathcal{F})/B^{s}_{k}(\mathcal{F})).  
\end{eqnarray*}
Also, from the dual complex $C^{\ast}(\mathcal{F}) := \mathrm{Hom}(C_{\ast}(\mathcal{F}), G)$, where $G = \Z$ or $\ZZ$, and $\delta_{\mathcal{F}}^{\ast}:=(\partial^{\mathcal{F}}_{\ast})^{\ast}$ we have the cohomology group of $\mathcal{F}$ by the standard homology theory:
\begin{eqnarray*}
 H^{\ast}(\mathcal{F}) & = & H_{f}^{\ast}(\mathcal{F})\oplus H_{s}^{\ast}(\mathcal{F}) \\
  & = & (Z_{f}^{\ast}(\mathcal{F})/B_{f}^{\ast}(\mathcal{F}))\oplus (Z_{s}^{\ast}(\mathcal{F})/B_{s}^{\ast}(\mathcal{F})).  
\end{eqnarray*}

\section{Duality theorem}
Let $(M, F), \mathcal{F}, \mathcal{M}$ and $\mathcal{R}$ be as in the previous section.  
\subsection{Statement of the duality theorem}
In this section, we prove the following theorem:

\vspace{2ex}\noindent
{\bf Theorem [Duality theorem]}
\textit{
 Let $\mathcal{F}^{\ast}$ be the backward foliation defined by the backward flow $F^{\ast}$.  With the $\ZZ$ coefficient, the following hold for $k=0, 1, \ldots , n-1$:
}
\begin{eqnarray*}
 H_{f}^{n-k-1}(\mathcal{F}^{\ast}) & \cong & H^{s}_{k}(\mathcal{F}) \\
 H_{s}^{n-k-1}(\mathcal{F}^{\ast}) & \cong & H^{f}_{k}(\mathcal{F}).  
\end{eqnarray*}
\textit{
In particular, we have
\[
 H^{n-k-1}(\mathcal{F}^{\ast})\cong H_{k}(\mathcal{F}).  
\]
Moreover, if there are no twisted compact leaves, i.e., all the associated round handles are trivial, then the same statements hold with $\Z$ coefficient.  
}

\subsection{The scheme of the proof}
In order to prove $H_{s}^{n-k-1}(\mathcal{F}^{\ast})\cong H^{f}_{k}(\mathcal{F})$, we will show $Z_{s}^{n-k-1}(\mathcal{F}^{\ast})\cong Z^{f}_{k}(\mathcal{F})$ and $B_{s}^{n-k-1}(\mathcal{F}^{\ast})\cong B^{f}_{k}(\mathcal{F})$.  Also, we will show $Z_{f}^{n-k-1}(\mathcal{F}^{\ast})\cong Z^{s}_{k}(\mathcal{F})$ and $B_{f}^{n-k-1}(\mathcal{F}^{\ast})\cong B^{s}_{k}(\mathcal{F})$ for $H_{f}^{n-k-1}(\mathcal{F}^{\ast})\cong H^{s}_{k}(\mathcal{F})$.  

First, we consider $H^{f}_{k}(\mathcal{F})$ and $H_{s}^{n-k-1}(\mathcal{F}^{\ast})$.  We have to compare the boundary homomorphism on $C^{f}_{k}(\mathcal{F})$ and the coboundary homomorphism on $C_{s}^{n-k-1}(\mathcal{F}^{\ast})$.  The targets of those homomorphisms are somewhat complicated, however, the boundary homomorphism is  \\
\hspace{10em}
$\partial^{f}_{k}:C^{f}_{k}(\mathcal{F})\rightarrow C^{f}_{k-1}(\mathcal{F})\oplus C^{s}_{k}(\mathcal{F})$
\\
and as it will be shown later, the coboundary homomorphism is \\
\hspace{10em}
$\delta_{s}^{n-k-1}:C_{s}^{n-k-1}(\mathcal{F}^{\ast})\rightarrow C_{s}^{n-k}(\mathcal{F}^{\ast})\oplus C_{f}^{n-k-1}(\mathcal{F}^{\ast})$.  
\\
Since those chain and cochain groups are free modules with natural bases, we have the representative matrices, thus we examine those entries.  We calculate the incident numbers and observe the difference of the matrices.  We find out that the difference is just a sign, i.e., the representative matrix of coboundary is just $\pm 1$ times of that of the boundary. Hence we obtain $Z^{f}_{k}(\mathcal{F})\cong Z_{s}^{n-k-1}(\mathcal{F}^{\ast})$.  

As for the boundary and the coboundary groups, we fix the indices of the targets of the boundary and the coboundary homomorphisms, 
\\
\hspace{8em}
$(\partial^{f}_{\ell +1}, \partial^{s}_{\ell +2}):C^{f}_{\ell +1}(\mathcal{F})\oplus C^{s}_{\ell +2}(\mathcal{F})\rightarrow C^{f}_{\ell}(\mathcal{F})$  
\\
and 
\\
\hspace{8em}
$(\delta_{s}^{n-\ell -2}, \delta_{f}^{n-\ell -3}):C_{s}^{n-\ell -2}(\mathcal{F}^{\ast})\oplus C_{f}^{n-\ell -3}(\mathcal{F}^{\ast})\rightarrow C_{s}^{n-\ell -1}(\mathcal{F}^{\ast})$.  
\\
By the argument similar to the case of cycles and cocycles, we have $B_{s}^{n-\ell -1}(\mathcal{F}^{\ast})\cong B^{f}_{\ell}(\mathcal{F})\quad (\ell =0, \ldots , n-1)$.  Therefore, we have $H_{s}^{n-k-1}(\mathcal{F}^{\ast})\cong H^{f}_{k}(\mathcal{F})$ for any $k=0, \ldots , n-1$.  

Also, we will show $Z_{f}^{n-\ell -1}(\mathcal{F}^{\ast})\cong Z^{s}_{\ell}(\mathcal{F})$ and $B_{f}^{n-\ell -1}(\mathcal{F}^{\ast})\cong B^{s}_{\ell}(\mathcal{F})$ and hence $H_{f}^{n-\ell -1}(\mathcal{F}^{\ast})\cong H^{s}_{k}(\mathcal{F})$ by an argument similar to the case of the fundamental chains.  

\subsection{Proof of the duality theorem}
Recall that the chain group $C_{\ast}(\mathcal{F})$ is a graded free module over $\Z$ or $\ZZ$ and those bases are bases of $C^{\mathcal{F}}_{\langle k+1\rangle}(\pi)=H_{k+1}(R(\pi), \partial_{-}R(\pi))$ and $C^{\mathcal{F}}_{\langle k\rangle}(\pi)=H_{k}(R(\pi), \partial_{-}R(\pi))$ for all $\pi\in P^{k}$ ($k=0, 1, \ldots, n-1$).  To compare the boundary and the coboundary homomorphisms, for any $\pi\in P^{k}$ and $\eta\in\mathrm{Blw}(\pi)$ we have to describe the component maps $\partial^{\mathcal{F}}_{\langle m\rangle}:C^{\mathcal{F}}_{\langle m\rangle}(\pi)\rightarrow C^{\mathcal{F}}_{\langle m-1\rangle}(\eta)$ and $\delta_{\mathcal{F}^{\ast}}^{\langle n-m-1\rangle}:C_{\mathcal{F}^{\ast}}^{\langle n-m-1\rangle}(\pi^{\ast})\rightarrow C_{\mathcal{F}^{\ast}}^{\langle n-m\rangle}(\eta^{\ast})$ ($m=k$ or $k+1$), where $\pi^{\ast}$ and $\eta^{\ast}$ denote the corresponding periodic orbits to $\pi$ and $\eta$ with respect to the backward flow.  The following lemma is clear.  
\begin{lem}
 For any $\pi\in P$ and any $m\in\Z \ (0\leq m\leq n-1)$, $C^{\mathcal{F}}_{\langle m\rangle}(\pi)=0$ if and only if $C_{\mathcal{F}^{\ast}}^{\langle n-m\rangle}(\pi^{\ast})=0$.  
\end{lem}
\noindent
By this lemma, it is enough to consider the case when the both of $C^{\mathcal{F}}_{\langle m\rangle}(\pi)$ and $C_{\mathcal{F}^{\ast}}^{\langle n-m\rangle}(\pi^{\ast})$ are non-trivial.  Suppose that $\mathrm{ind}(\pi)=k$ and $\mathrm{ind}(\eta)=\ell$, and we write $R^{k}=R(\pi), R^{\ell}=R(\eta)$.  Now we will fix a generator of $C^{\mathcal{F}}_{\langle m\rangle}(\pi)=H_{m}(R^{k}, \partial_{-}R^{k})\cong\mathscr{H}^{\mathcal{R}}_{m}(\pi)$\ $(m =k$ or $k+1)$.  Recall that round $k$-handle is the total space of a Whitney sum of two disk bundles $E^{k}$ and $E^{n-k-1}$ over $S^{1}$, where $E^{p}$ denotes a $D^{p}$-bundle.  We write $R^{k}=E^{k}\oplus E^{n-k-1}$ and $\pi_{E^{k}\oplus E^{n-k-1}}:E^{k}\oplus E^{n-k-1}\rightarrow S^{1}$.  First, we define subsets $L_{\pm}(R^{k})\subset K_{\pm}(R^{k})\subset R^{k}$ as follows:
\begin{eqnarray*}
 K_{-}(R^{k}) & := & E^{k}\oplus 0 \\
 L_{-}(R^{k}) & := & (\pi_{E^{k}\oplus 0})^{-1}(1) \\
 K_{+}(R^{k}) & := & 0\oplus E^{n-k-1} \\
 L_{+}(R^{k}) & := & (\pi_{0\oplus E^{n-k-1}})^{-1}(1) 
\end{eqnarray*}
where we consider $S^{1}$ as $\partial D^{2}\subset D^{2}\subset\C$.  From now on, let $G$ be $\Z$ or $\ZZ$, depending on the situation that the round handles appearing are twisted or not:
\[
 G = \left\{\begin{array}{cl}
      \Z & \text{if the RHD is simple, } \\
      \ZZ & \text{otherwise,}
	    \end{array} \right.  
\]
and we assume that the coefficient group of the homology groups are $G$.  
The following lemma can be seen by the standard homology theory:  
\begin{lem}\label{Intersection}
 \begin{enumerate}
  \item The homology class $[K_{-}(R^{k}), \partial K_{-}(R^{k})]$ generates $C^{\mathcal{F}}_{\langle k+1\rangle}(\pi)\\ 
=H_{k+1}(R^{k}, \partial_{-}R^{k})\cong G$. 
  \item The homology class $[L_{-}(R^{k}), \partial L_{-}(R^{k})]$ generates $C^{\mathcal{F}}_{\langle k\rangle}(\pi)=H_{k}(R^{k}, \partial_{-}R^{k})\cong G$. 
  \item\label{I_f} For any $z\in H_{k+1}(R^{k}, \partial_{-}R^{k})$, $z=\ell [K_{-}(R^{k}), \partial K_{-}(R^{k})]\ (\ell\in G)$ if and only if $\ell = I_{f}^{k}(z, [L_{+}(R^{k}), \partial L_{+}(R^{k})])$, where $I_{f}^{k}:H_{k+1}(R^{k}, \partial_{-}R^{k})\times H_{n-k-1}(R^{k}, \partial_{+}R^{k})\rightarrow G$ is the intersection form.  
  \item\label{I_s} For any $w\in H_{k}(R^{k}, \partial_{-}R^{k})$, $w=m [L_{-}(R^{k}), \partial L_{-}(R^{k})]\ (m\in G)$ if and only if $m= I_{s}^{k}(w, [K_{+}(R^{k}), \partial K_{+}(R^{k})])$, where $I_{s}^{k}:H_{k}(R^{k}, \partial_{-}R^{k})\times H_{n-k}(R^{k}, \partial_{+}R^{k})\rightarrow G$ is the intersection form.  
 \end{enumerate}
\end{lem}
For the simplicity, we write
\[
 \begin{array}{lcl}
  \kappa_{\pm}^{k} & = & [K_{\pm}(R^{k}), \partial K_{\pm}(R^{k})] \\
  \kappa_{\pm}^{k\ast} & = & [K_{\pm}(R^{k\ast}), \partial K_{\pm}(R^{k\ast})] \\
  \lambda_{\pm}^{k} & = & [L_{\pm}(R^{k}), \partial L_{\pm}(R^{k})] \\
  \lambda_{\pm}^{k\ast} & = & [L_{\pm}(R^{k\ast}), \partial L_{\pm}(R^{k\ast})]   
 \end{array}
\]
where $R^{k\ast}$ denotes the round handle $R^{k}$ with the backward flow.  

\begin{remark}
The index of the backward round handle $R^{k\ast}$ is equal to $n-k-1$ and $R^{k}=R^{k\ast}=R^{n-k-1}$ as sets.  Note that the orientations of $R^{k}$ and $R^{k\ast}$ may be changed.  Thus we have
\begin{eqnarray*}
 K_{-}(R^{k\ast}) & = & K_{+}(R^{k}), \\
 L_{-}(R^{k\ast}) & = & L_{+}(R^{k}), \\
 K_{+}(R^{k\ast}) & = & K_{-}(R^{k}), \\
 L_{+}(R^{k\ast}) & = & L_{-}(R^{k}).  
\end{eqnarray*}
\end{remark}
Now, we proceed to describe and to compare the boundary and the coboundary homomorphisms in the cases divided by the homology dimensions $\langle k\rangle \text{ or }\langle k+1\rangle$ of $R^{k}=R(\pi)\ $ and the indices of $R^{\ell}=R(\eta)\ (\ell =\mathrm{ind}(\eta)=k-2, k-1\text{ or } k)$.  From now on, we fix a periodic orbit $\pi$ of index $k$ as the source of the boundary homomorphism.  

\begin{case}\label{f2f}
 From the fundamental $k$-chain to a fundamental $(k-1)$-chain.  
\end{case}

We consider the component of the boundary $\partial^{f}_{k}:H_{k+1}(R^{k}, \partial_{-}R^{k})\rightarrow H_{k}(R^{k-1}, \partial_{-}R^{k-1})$.  Let $\partial^{f}_{k}(\kappa_{-}^{k})=[R_{f}^{k};R_{f}^{k-1}]\kappa_{-}^{k-1}$, where $[R_{f}^{k};R_{f}^{k-1}]\in G$ which is called the {\em incident number} between the fundamental chain of $R^{k}$ and the fundamental chain of $R^{k-1}$.  Note that we add the subscripts $f$ or $s$, here,  to the symbols of the round handles for the clarity.  
Next we consider the backward flow from $R^{k-1\ast}$ to $R^{k\ast}$.  In the backward flow,  $R^{k-1\ast}$ is a round $(n-k)$-handle and $R^{k\ast}$ is a round $(n-k-1)$-handle, and the boundary is $\partial^{s}_{n-k}:H_{n-k}(R^{k-1\ast}, \partial_{-}R^{k-1\ast})\rightarrow H_{n-k-1}(R^{k\ast}, \partial_{-}R^{k\ast})$.  Suppose that $\partial^{s}_{n-k}(\lambda_{-}^{k-1\ast})=[R_{s}^{k-1\ast};R_{s}^{k\ast}]\lambda_{-}^{k\ast}$.  

\begin{lem}\label{Case1_lemma}
 The above two incident numbers are equal up to sign. In fact, $[R_{s}^{k-1\ast};R_{s}^{k\ast}]=(-1)^{k+1}[R_{f}^{k};R_{f}^{k-1}]$ if the RHD is simple and $[R_{s}^{k-1\ast};R_{s}^{k\ast}]=[R_{f}^{k};R_{f}^{k-1}]$ otherwise.  
\end{lem}

\begin{proof}
By Lemma \ref{Intersection} (\ref{I_f}) and (\ref{I_s}), we have $[R^{k}_{f};R^{k-1}_{f}]=I_{f}^{k-1}(\partial^{f}_{k}(\kappa_{-}^{k}), \lambda_{+}^{k-1})$ and $[R^{k-1\ast}_{s};R^{k\ast}_{s}]=I_{s}^{k\ast}(\partial^{s}_{n-k}(\lambda_{-}^{k-1\ast}), \kappa_{+}^{k\ast})$, where $I_{f}^{k-1}:H_{k}(R^{k-1}, \partial_{-}R^{k-1})\times H_{n-k}(R^{k-1}, \partial_{+}R^{k-1})\rightarrow G$ and $I_{s}^{k\ast}:H_{n-k-1}(R^{k\ast}, \partial_{-}R^{k\ast})\times H_{k+1}(R^{k\ast}, \partial_{+}R^{k\ast})\rightarrow G$ are intersection forms.  As cycles which represent $\partial^{f}_{k}(\kappa_{-}^{k})$ and $\partial^{s}_{n-k}(\lambda_{-}^{k-1\ast})$, we may consider $(\partial K_{-}(R^{k})\cdot\epsilon)\cap(R^{k-1}, \partial_{-}R^{k-1})$ and \\ 
$(\partial L_{-}(R^{k-1\ast})\cdot\epsilon)\cap (R^{k\ast}, \partial_{-}R^{k\ast})$ for a small enough $\epsilon >0$.  Here, we note that $x\cdot\epsilon$ in the backward flow implies $F(x, -\epsilon)$.  By the transversality, $(L_{-}(R^{k-1\ast})\cdot\epsilon)\cap (K_{+}(R^{k\ast})\cdot(-\epsilon))$ is a union of simple arcs each of which connects a point in $(\partial L_{-}(R^{k-1\ast})\cdot\epsilon)\cap K_{+}(R^{k\ast})$ to a point in $L_{-}(R^{k-1\ast})\cap(\partial K_{+}(R^{k\ast})\cdot (-\epsilon))$.  Thus the points of $(\partial L_{-}(R^{k-1\ast})\cdot\epsilon)\cap K_{+}(R^{k\ast})$ and $L_{-}(R^{k-1\ast})\cap(\partial K_{+}(R^{k\ast})\cdot (-\epsilon))=L_{+}(R^{k-1})\cap (\partial K_{-}(R^{k})\cdot\epsilon )$ correspond to each other bijectively.  In the case of non-simple RHD, this shows the assertion of the lemma.  

In the following, we only consider the case of simple RHD.  
Then, counting the number of sighed intersection points in the cycle level, we have the following (the second equality is shown by Sublemma \ref{intchange_del} below):
\begin{eqnarray*}
 I_{s}^{k\ast}(\partial^{s}_{n-k}(\lambda_{-}^{k-1\ast}),\ \kappa_{+}^{k\ast}) 
  & = & \langle\partial L_{-}(R^{k-1\ast})\cdot\epsilon,\ K_{+}(R^{k\ast})\rangle_{R^{k\ast}} \\
  & = & (-1)^{k+1}\langle L_{-}(R^{k-1\ast}),\ \partial K_{+}(R^{k\ast})\cdot (-\epsilon)\rangle_{R^{k-1\ast}} \\
  & = & (-1)^{k+1}\langle L_{+}(R^{k-1}),\ \partial K_{-}(R^{k})\cdot(-\epsilon)\rangle_{R^{k-1\ast}} \\
  & = & (-1)^{k+1}\langle  \partial K_{-}(R^{k})\cdot\epsilon,\ L_{+}(R^{k-1})\rangle_{R^{k-1}} \\
  & = & (-1)^{k+1}I_{f}^{k-1}(\partial^{f}_{k}(\kappa_{-}^{k}),\ \lambda_{+}^{k-1})
\end{eqnarray*}
where, e.g., $\langle Z_{1},\ Z_{2}\rangle_{X}$ denotes the signed number of intersection points of $\ell$-cycle $Z_{1}$ and $(n-\ell)$-cycle $Z_{2}$ in general position in $X$.  
\begin{sublem}\label{intchange_del}
 $\langle\partial L_{-}(R^{k-1\ast})\cdot\epsilon,\ K_{+}(R^{k\ast})\rangle_{R^{k\ast}}
 = (-1)^{k+1}\langle L_{-}(R^{k-1\ast}),\ \partial K_{+}(R^{k\ast})\cdot (-\epsilon)\rangle_{R^{k-1\ast}}$
\end{sublem}
\begin{proof}

Let $e$ denote a component of $(L_{-}(R^{k-1\ast})\cdot\epsilon)\cap (K_{+}(R^{k\ast})\cdot(-\epsilon))$ and $\partial e=\{ a, b\}$ where $a\in\partial L_{-}(R^{k-1\ast})\cdot\epsilon$ and $b\in\partial K_{+}(R^{k\ast})\cdot (-\epsilon)$.  Let $c$ denote the midpoint $e\cap\partial_{-}R^{k-1\ast}=e\cap\partial_{+}R^{K\ast}$  We take a frame $[v_{1}v_{2}\cdots v_{n-k}]$ of $L_{-}(R^{k-1\ast})$ near $e$ in $R^{k-1\ast}$ such that $v_{1},\ldots , v_{n-k-1}$ are tangent to $\partial L_{-}(R^{k-1\ast})$ and $v_{n-k}$ points inward tangent to $e$ at $c$.  Moreover, let $[\bar{v}_{1}\cdots\bar{v}_{n-k-1}]$ denote the frame of $\partial L_{-}(R^{k-1\ast})\cdot\epsilon$ at $a$ in $R^{k\ast}$ induced from $[v_{1}\cdots v_{n-k-1}]$ by the flow.  Also we take a frame $[w_{0}w_{1}\cdots w_{k}]$ of $K_{+}(R^{k\ast})$ near $e$ in $R^{k\ast}$ where $w_{0}$ points parallel to the zero-section (=$S^{1}$-factor) and $w_{1},\ldots, w_{k-1}$ are tangent to $\partial K_{+}(R^{k\ast})$ and $w_{k}$ points inward tangent to $e$ at $c$.  Similarly as above, let $[\bar{w}_{0}\cdots\bar{w}_{k-1}]$ denote the frame of $\partial K_{+}(R^{k\ast})\cdot(-\epsilon)$ at $b$ in $R^{k-1\ast}$ induced from $[w_{0}\cdots w_{k-1}]$. 
Then we have the frame $[\bar{v}_{1}\cdots\bar{v}_{n-k-1}][w_{0}w_{1}\cdots w_{k}]$ which determines a local orientation at $a$.  On the other hand, the frame $[v_{1}\cdots v_{n-k}][\bar{w}_{0}\bar{w}_{1}\cdots\bar{w}_{k-1}]$ determines another local orientation at $b$.  At $\epsilon =0$, $v_{n-k}$ and $w_{k}$ can be taken as $w_{k}=-v_{n-k}$ and we have two oriented frames at $c$, $[v_{1}\cdots v_{n-k-1}][w_{0}w_{1}\cdots w_{k}]$ and $[v_{1}\cdots v_{n-k}][w_{0}w_{1}\cdots w_{k-1}]$.  Therefore, rearranging $w_{k}$ of $[v_{1}\cdots v_{n-k-1}w_{0}w_{1}\cdots w_{k}]$ to the $(n-k)$th position, we have the frame $[v_{1}\cdots v_{n-k-1}w_{k}w_{0}\cdots w_{k-1}]$ which is the opposite frame of $[v_{1}\cdots v_{n-k}w_{0}\cdots w_{k-1}]$ at $c$, i.e., on the boundary $\partial_{-}R^{k-1\ast}\cap\partial_{+}R^{k\ast}$.  Hence the signs differ by the factor $(-1)^{k+1}$.  
\end{proof}
\noindent
Thus, we have $[R_{f}^{k};R_{f}^{k-1}]=(-1)^{k+1}[R_{s}^{k-1\ast};R_{s}^{k\ast}]$.  
\end{proof}
Hence the following diagram is commutative 
\begin{equation*}
 \vcenter{
  \xymatrix@C=30pt{
   C^{\mathcal{F}}_{\langle k+1\rangle}(\pi) \ar[r]^{\partial^{\mathcal{F}}}
 & C^{\mathcal{F}}_{\langle k\rangle}(\eta) \\
   G \ar[u]^{\cong} & G \ar[u]_{\cong} \\
   C^{\mathcal{F^{\ast}}}_{\langle n-k\rangle}(\eta^{\ast}) \ar[r] \ar[u]^{\cong} 
 & C^{\mathcal{F^{\ast}}}_{\langle n-k-1\rangle}(\pi^{\ast}) \ar[u]_{\cong} 
}
}
\end{equation*}
where the horizontal arrow below is $(-1)^{k+1}\partial^{\mathcal{F^{\ast}}}$ or $\partial^{\mathcal{F^{\ast}}}$, depending on the types of RHD.  This implies that the dual homomorphism $\delta_{\mathcal{F}^{\ast}}^{\langle n-k-1\rangle}=(\partial^{\mathcal{F^{\ast}}}_{\langle n-k\rangle})^{\ast}:C_{\mathcal{F^{\ast}}}^{\langle n-k-1\rangle}(\pi^{\ast})\rightarrow C_{\mathcal{F^{\ast}}}^{\langle n-k\rangle}(\eta^{\ast})$ is also the $(-1)^{k+1}[R_{f}^{k};R_{f}^{k-1}]$ or $[R_{f}^{k};R_{f}^{k-1}]$ times map from $G$ to itself with respect to the natural farmings.  This shows the comparison in the Case \ref{f2f}.  

In the following, we omit the argument in the case of non-simple RHD since it is the same as in Case 1.  

\begin{case}\label{f2s}
 From the fundamental $k$-chain to a sectional $k$-chain.  
\end{case}

This is the case of $\partial^{f}_{k}:H_{k+1}(R_{1}^{k}, \partial_{-}R_{1}^{k})\rightarrow H_{k}(R_{2}^{k}, \partial_{-}R_{2}^{k})$.  As in the argument in the Case \ref{f2f}, suppose that $\partial^{f}_{k}(\kappa_{1-}^{k})=[R_{1f}^{k};R_{2s}^{k}]\lambda_{2-}^{k}$ where $\kappa_{i\pm}^{k}, \lambda_{i\pm}^{k}$ denote the homology classes corresponding to $R_{i}^{k}\ (i=1, 2)$.  As for the backward flow, the component of the boundary homomorphism is $\partial^{f}_{n-k-1}:H_{n-k}(R_{2}^{k\ast}, \partial_{-}R_{2}^{k\ast})\rightarrow H_{n-k-1}(R_{1}^{k\ast}, \partial_{-}R_{1}^{k\ast})$ and suppose $\partial^{f}_{n-k-1}(\kappa_{2-}^{k\ast})=[R_{2f}^{k\ast};R_{1s}^{k\ast}]\lambda_{1-}^{k\ast}$.  

\begin{lem}
 $[R_{2f}^{k\ast};R_{1s}^{k\ast}]=(-1)^{k+1}[R_{1f}^{k};R_{2s}^{k}]$.  
\end{lem}

\begin{proof}[Sketch of the froof]
 We proceed as in the proof of Lemma \ref{Case1_lemma}.  We have $[R_{1f}^{k};R_{2s}^{k}]=I^{k}_{s}(\partial^{f}_{k}(\kappa_{1-}^{k}), \kappa_{2+}^{k})$ and $[R_{2f}^{k\ast};R_{1s}^{k\ast}]=I^{k\ast}_{s}(\partial^{f}_{n-k-1}(\kappa_{2-}^{k\ast}), \kappa_{1+}^{k\ast})$.  Then, similarly as in the proof of Sublemma \ref{intchange_del}, we have the following: 
\begin{eqnarray*}
 I_{s}^{k\ast}(\partial^{f}_{n-k-1}(\kappa_{2-}^{k\ast}),\ \kappa_{1+}^{k\ast})
  & = & \langle\partial K_{-}(R_{2}^{k\ast})\cdot\epsilon,\ K_{+}(R_{1}^{k\ast})\rangle_{R_{1}^{k\ast}} \\
  & = & (-1)^{k+1}\langle K_{-}(R_{2}^{k\ast}),\ \partial K_{+}(R_{1}^{k\ast})\cdot (-\epsilon)\rangle_{R_{2}^{k\ast}} \\
  & = & (-1)^{k+1}\langle K_{+}(R_{2}^{k}),\ \partial K_{-}(R_{1}^{k})\cdot (-\epsilon)\rangle_{R_{2}^{k\ast}} \\
  & = & (-1)^{k+1}\langle \partial K_{-}(R_{1}^{k})\cdot\epsilon,\ K_{+}(R_{2}^{k})\rangle_{R_{2}^{k}} \\
  & = & (-1)^{k+1}I_{s}^{k}(\partial^{f}_{k}(\kappa_{1-}^{k}), \kappa_{2+}^{k})
\end{eqnarray*}
\end{proof}
The cochain map $\delta_{s}^{n-k-1}=(\partial^{f}_{n-k-1})^{\ast}:(C^{\mathcal{F}}_{\langle n-k-1\rangle}(\pi^{k\ast}))^{\ast}\rightarrow (C^{\mathcal{F}}_{\langle n-k\rangle}(\eta^{k\ast}))^{\ast}$ is also $[R_{2f}^{k\ast};R_{1s}^{k\ast}]$ times map from $G$ to itself.  This shows the comparison in the Case 2. 

The rest of the comparisons are two cases, each of which we only show a sketch of the proof since the proofs proceed similarly as in the Case \ref{f2f}, \ref{f2s}.  

\begin{case}\label{s2s}
 From the sectional $k$-chain to a sectional $(k-1)$-chain.  
\end{case}

This is the case of $\partial^{s}_{k}:H_{k}(R^{k}, \partial_{-}R^{k})\rightarrow H_{k-1}(R^{k-1}, \partial_{-}R^{k-1})$ and suppose $\partial^{s}_{k}(\lambda_{-}^{k})=[R_{s}^{k};R_{s}^{k-1}]\lambda_{-}^{k-1}$.  For the backward flow, the boundary is $\partial^{f}_{n-k}:H_{n-k+1}(R^{k-1\ast}, \partial_{-}R^{k-1\ast})\rightarrow H_{n-k}(R^{k\ast}, \partial_{-}R^{k\ast})$ and suppose $\partial^{f}_{n-k}(\kappa_{-}^{k-1\ast})=[R_{f}^{k-1\ast};R_{f}^{k\ast}]\kappa_{-}^{k\ast}$.  

\begin{lem}
 $[R_{f}^{k-1\ast};R_{f}^{k\ast}]=(-1)^{k}[R_{s}^{k};R_{s}^{k-1}]$.
\end{lem}

\begin{proof}[Sketch of the proof]
 \begin{eqnarray*}
  [R_{f}^{k-1\ast};R_{f}^{k\ast}]
   & = & I_{f}^{k\ast}(\partial^{f}_{n-k}(\kappa_{-}^{k-1\ast}),\ \lambda_{+}^{k\ast}) \\
   & = & \langle\partial K_{-}(R^{k-1\ast})\cdot\epsilon,\ L_{+}(R^{k\ast})\rangle_{R^{k\ast}} \\
   & = & (-1)^{k}\langle K_{-}(R^{k-1\ast}),\ \partial L_{+}(R^{k\ast})\cdot (-\epsilon)\rangle_{R^{k-1\ast}}  \\
   & = & (-1)^{k}\langle K_{+}(R^{k-1}),\ \partial L_{-}(R^{k})\cdot(-\epsilon)\rangle_{R^{k-1\ast}} \\
   & = & (-1)^{k}\langle  \partial L_{-}(R^{k})\cdot\epsilon,\ K_{+}(R^{k-1})\rangle_{R^{k-1}} \\
   & = & (-1)^{k}I_{s}^{k-1}(\partial^{s}_{k}(\lambda_{-}^{k}),\ \kappa_{+}^{k-1}) \\
   & = & (-1)^{k}[R_{s}^{k};R_{s}^{k-1}]
 \end{eqnarray*}
\end{proof}

Thus the difference factor between $\partial^{\mathcal{F}}_{\langle k\rangle}:C^{\mathcal{F}}_{\langle k\rangle}(\pi)\rightarrow C^{\mathcal{F}}_{\langle k-1\rangle}(\eta)$ and $\delta_{\mathcal{F}^{\ast}}^{\langle n-k\rangle}:C_{\mathcal{F}^{\ast}}^{\langle n-k\rangle}(\pi^{\ast})\rightarrow C_{\mathcal{F}^{\ast}}^{\langle n-k+1\rangle}(\eta^{\ast})$ is $(-1)^{k}$.  

\begin{case}\label{s2f}
 From the sectional $k$-chain to a fundamental $(k-2)$-chain.  
\end{case}

In this case the boundary homomorphisms are $\partial^{s}_{k}:H_{k}(R^{k}, \partial_{-}R^{k})\rightarrow H_{k-1}(R^{k-2}, \partial_{-}R^{k-2})$ and $\partial^{s}_{n-k+1}:H_{n-k+1}(R^{k-2\ast}, \partial_{-}R^{k-2\ast})\rightarrow H_{n-k}(R^{k\ast}, \partial_{-}R^{k\ast})$.  Suppose $\partial^{s}_{k}(\lambda_{-}^{k})=[R_{s}^{k};R_{f}^{k-2}]\kappa_{-}^{k-2}$ and 
$\partial^{s}_{n-k+1}(\lambda_{-}^{k-2\ast})=[R_{s}^{k-2\ast};R_{f}^{k\ast}]\kappa_{-}^{k\ast}$.  

\begin{lem}
 $[R_{s}^{k-2\ast};R_{f}^{k\ast}]=(-1)^{k}[R_{s}^{k};R_{f}^{k-2}]$. 
\end{lem} 

\begin{proof}[Sketch of the proof]
\begin{eqnarray*}
 [R_{s}^{k-2\ast};R_{f}^{k\ast}]
  & = & I_{f}^{k\ast}(\partial^{s}_{n-k+1}(\lambda_{-}^{k-2\ast}),\ \lambda_{+}^{k\ast}) \\
  & = & \langle\partial L_{-}(R^{k-2\ast})\cdot\epsilon,\ L_{+}(R^{k\ast})\rangle_{R^{k\ast}} \\
  & = & (-1)^{k}\langle L_{-}(R^{k-2\ast}),\ \partial L_{+}(R^{k\ast})\cdot (-\epsilon)\rangle_{R^{k-2\ast}} \\
  & = & (-1)^{k}\langle \partial L_{-}(R^{k})\cdot\epsilon,\ L_{+}(R^{k-2})\rangle_{R^{k-2}} \\
  & = & (-1)^{k}I_{f}^{k-2}(\partial^{s}_{k}(\lambda_{-}^{k}),\ \lambda_{+}^{k-2}) \\
  & = & (-1)^{k}[R_{s}^{k};R_{f}^{k-2}].  
\end{eqnarray*}
\end{proof}
Thus the difference factor between $\partial^{\mathcal{F}}_{\langle k\rangle}:C^{\mathcal{F}}_{\langle k\rangle}(\pi)\rightarrow C^{\mathcal{F}}_{\langle k-1\rangle}(\eta)$ and $\delta_{\mathcal{F}^{\ast}}^{\langle n-k\rangle}:C_{\mathcal{F}^{\ast}}^{\langle n-k\rangle}(\pi^{\ast})\rightarrow C_{\mathcal{F}^{\ast}}^{\langle n-k+1\rangle}(\eta^{\ast})$ is $(-1)^{k}$.  

\noindent
This completes the comparisons of the component maps of the boundary and the coboundary homomorphisms.  
\vspace{2ex}

Summarising the results, we have the following list.  
\[
\begin{array}{l|l|l|c}
  & \text{h-dim}(\pi) & \mathrm{ind}(\eta) & \text{exponent of difference factor} \\ \hline
 \text{Case 1} & \langle k+1\rangle\ \mathrm{(fund.)} & k-1\ \mathrm{(fund.)} & k+1 \\
 \text{Case 2} & \langle k+1\rangle\ \mathrm{(fund.)} & k\ \mathrm{(sect.)} & k+1 \\
 \text{Case 3} & \langle k\rangle\ \mathrm{(sect.)} & k-1\ \mathrm{(sect.)} & k \\
 \text{Case 4} & \langle k\rangle\ \mathrm{(sect.)} & k-2\ \mathrm{(fund.)} & k 
\end{array} 
\]
\noindent
Note that at the fundamental $k$-chains both of the difference factors in the Case \ref{f2f} and \ref{f2s} are $(-1)^{k+1}$, and at the sectional $k$-chains they are $(-1)^{k}$.  Hence with respect to the natural framing the representative matrix of $\delta_{s}^{n-k-1}:C_{s}^{n-k-1}(\mathcal{F}^{\ast})\rightarrow C_{s}^{n-k}(\mathcal{F}^{\ast})\oplus C_{f}^{n-k-1}(\mathcal{F}^{\ast})$ is the $(-1)^{k+1}$ times of that of $\partial^{f}_{k}:C^{f}_{k}(\mathcal{F})\rightarrow C^{f}_{k-1}(\mathcal{F})\oplus C^{s}_{k}(\mathcal{F})$.  Therefore we have $Z_{s}^{n-k-1}(\mathcal{F}^{\ast})\cong Z^{f}_{k}(\mathcal{F})$.  A similar argument on $\delta_{f}^{n-k-1}:C_{f}^{n-k-1}(\mathcal{F}^{\ast})\rightarrow C_{f}^{n-k}(\mathcal{F}^{\ast})\oplus C_{s}^{n-k+1}(\mathcal{F}^{\ast})$ and $\partial^{s}_{k}:C^{s}_{k}(\mathcal{F})\rightarrow C^{s}_{k-1}(\mathcal{F})\oplus C^{f}_{k-2}(\mathcal{F})$ shows $Z_{f}^{n-k-1}(\mathcal{F}^{\ast})\cong Z^{s}_{k}(\mathcal{F})$.  

Regarding the boundaries and the coboundaries, we need the following list in which the index $\ell$ of the target $\eta$ is fixed.
\[
\begin{array}{l|l|l|c}
  & \mathrm{ind}(\pi) & \text{h-dim}(\eta) & \text{exponent of difference factor} \\ \hline
 \text{Case 1} & \ell +1\ \mathrm{(fund.)} & \langle \ell +1\rangle\ \mathrm{(fund.)} &  \ell +2 \\
 \text{Case 2} & \ell\ \mathrm{(fund.)} & \langle \ell\rangle\ \mathrm{(sect.)} & \ell +1 \\
 \text{Case 3} & \ell +1\ \mathrm{(sect.)} & \langle \ell\rangle\ \mathrm{(sect.)} & \ell +1 \\
 \text{Case 4} & \ell +2\ \mathrm{(sect.)} & \langle \ell +1\rangle\ \mathrm{(fund.)} & \ell +2
\end{array} 
\]
\noindent
Since the difference factors of the boundary homomorphisms into fundamental chains of $\eta$ are equal, the representative matrices of $(\partial^{f}_{\ell +1},\ \partial^{s}_{\ell +2}):C^{f}_{\ell +1}(\mathcal{F})\oplus C^{s}_{\ell +2}(\mathcal{F})\rightarrow C^{f}_{\ell}(\mathcal{F})$ and $(\delta_{s}^{n-\ell -2},\ \delta_{f}^{n-\ell -3}):C_{s}^{n-\ell -2}(\mathcal{F}^{\ast})\oplus C_{f}^{n-\ell -3}(\mathcal{F}^{\ast})\rightarrow C_{s}^{n-\ell -1}(\mathcal{F}^{\ast})$ are equal up to sign.  Thus, it follows $B_{s}^{n-\ell -1}(\mathcal{F}^{\ast})\cong B^{f}_{\ell}(\mathcal{F})$.  Hence we have $H_{s}^{n-k-1}(\mathcal{F}^{\ast})\cong H^{f}_{k}(\mathcal{F})$ for any $k=0, \ldots , n-1$.  

At the sectional chains, we have the similar situation.  Now by an argument similar to the case of the fundamental chains, it is easy to see $Z_{f}^{n-k-1}(\mathcal{F}^{\ast})\cong Z^{s}_{k}(\mathcal{F})$ and $B_{f}^{n-\ell -1}(\mathcal{F}^{\ast})\cong B^{s}_{\ell}(\mathcal{F})$.  Hence we have $H_{f}^{n-k-1}(\mathcal{F}^{\ast})\cong H^{s}_{k}(\mathcal{F})$ for any $k=0, \ldots , n-1$.  This completes the proof of the duality theorem.

\section*{acknowledgements}
The author would like to express gratitude to Atsushi Sato for his support.  This research could not be completed without his crucial comments and questions.  Also, the author would like to thank the Department of Mathematics of Chuo University for the kind permission for him to use the facilities.

\end{document}